\documentclass[opre,nonblindrev]{informs3}

\SingleSpacedXI

\usepackage{caption}
 \usepackage{subcaption}
\usepackage{multicol}
\usepackage{multirow}
\usepackage{hyperref}
\usepackage{bm}

\usepackage{endnotes}
\let\footnote=\endnote

\usepackage{natbib}
 \bibpunct[, ]{(}{)}{,}{a}{}{,}%
 \def\bibfont{\small}%
 \def\bibsep{\smallskipamount}%
 \def\bibhang{24pt}%
\TheoremsNumberedThrough     
\ECRepeatTheorems

\EquationsNumberedThrough    

\definecolor{bluecolor}{rgb}{0,0,1}

\definecolor{redcolor}{rgb}{1,0,0}

\begin{document}



\RUNTITLE{Fabrication-Adaptive Optimization}

\TITLE{Fabrication-Adaptive Optimization, with an Application to Photonic Crystal Design}

\ARTICLEAUTHORS{
\AUTHOR{Han Men}
\AFF{Department of Aeronautics \& Astronautics, Massachusetts Institute of Technology, Cambridge, MA 02139, \EMAIL{abbymen@mit.edu}} 
\AUTHOR{Robert M. Freund}
\AFF{Sloan School of Management, Massachusetts Institute of Technology, Cambridge, MA 02139, \EMAIL{rfreund@mit.edu}}
\AUTHOR{Ngoc C. Nguyen, Joel Saa-Seoane, Jaime Peraire}
\AFF{Department of Aeronautics \& Astronautics, Massachusetts Institute of Technology, Cambridge, MA 02139, \EMAIL{cuongng@mit.edu}, \EMAIL{jsaa@mit.edu}, \EMAIL{peraire@mit.edu}} 
} 

\ABSTRACT{It is often the case that the computed optimal solution of an optimization problem cannot be implemented directly, irrespective of data accuracy, due to either ({\em i}) technological limitations (such as physical tolerances of machines or processes), ({\em ii}) the deliberate simplification of a model to keep it tractable (by ignoring certain types of constraints that pose computational difficulties), and/or ({\em iii}) human factors (getting people to ``do'' the optimal solution).  Motivated by this observation, we present a modeling paradigm called ``fabrication-adaptive optimization'' for treating issues of implementation/fabrication.  We develop computationally-focused theory and algorithms, and we present computational results for incorporating considerations of implementation/fabrication into constrained optimization problems that arise in photonic crystal design.  The fabrication-adaptive optimization framework stems from the robust regularization of a function.  When the feasible region is not a normed space (as typically encountered in application settings), the fabrication-adaptive optimization framework typically yields a non-convex optimization problem.  (In the special case where the feasible region is a finite-dimensional normed space, we show that fabrication-adaptive optimization can be re-cast as an instance of modern robust optimization.)  We study a variety of problems with special structures on functions, feasible regions, and norms, for which computation is tractable, and develop an algorithmic scheme for solving these problems in spite of the challenges of non-convexity.  We apply our methodology to compute fabrication-adaptive designs of two-dimensional photonic crystals with a variety of prescribed features.}



\maketitle

%


\section{Introduction:  Problem Statement, Preliminaries, and Computational Aspirations}
Consider a general constrained optimization problem of the form:
\begin{equation}\label{eq_nominal}
\begin{array}{ll}
z^{\ast} = & \min\limits_{x} f(x)   \\
  &  \mbox{s.t.} \ x\in S,
\end{array}
\end{equation}
where $S\subseteq \mathbb{R}^n$ is the feasible region, and $f(\cdot): \mathbb{R}^n \rightarrow \mathbb{R}$ is the objective to be optimized.  The context of (\ref{eq_nominal}) may be as diverse as portfolio optimization \citep{port} where $x_i$ is the number of shares to be invested in asset $i$, to optimal microstructure material design where $x_i$ is the concentration of a dialectric material in pixel (or voxel) $i$ of  a discretized physical region as in \cite{men2010bandgap,men2011design}.  Let $x^{\ast}$ be an optimal solution of (\ref{eq_nominal}).  The notion of \emph{fabrication-adaptivity}, or perhaps more generally \emph{implementation-adaptivity}, has to do with the concern that while the data and other descriptors of the problem may be quite accurate, it may be generically implausible to implement the optimal solution $x^{\ast}$ exactly.  Some reasons for this may include:

\begin{enumerate}
\item Technological limitations.  For example, the production or fabrication technology might not be able to fabricate the product exactly according to the plan specified in ${x}^{\ast}$, perhaps due to limitations of machine tolerances,
\item Deliberate simplifications.  The model \eqref{eq_nominal} may be a deliberate simplification of the real problem in order for the optimization model to be computationally tractable.  For example, it may be computationally prohibitive to include odd-lot constraints in a portfolio optimization model, or connectivity constraints in a microstructure material design model, etc., or
\item Human factors.  It may be implausible to assume that people will ``do'' $x^{\ast}$ precisely as the optimization model prescribes.
\end{enumerate}\medskip

Indeed, the application that has given rise to this line of study is a nonlinear optimization problem arising in microstructure material design \citep{men2010bandgap}, where $S:=[x_{\min}, x_{\max}]^n = \{{x}\in\mathbb{R}^n: x_{\min} e \le {x} \le x_{\max} {e}\}$ is a hypercube with $n>>0$, and the component values $x_i$ represent the permittivity of a dielectric material at pixel (or voxel) $i$ for $i = 1,\ldots,n$.  The values $x_{\min}$ and $x_{\max}$ correspond to permittivity constants for air and the dialectric material (e.g., gallium arsenide), respectively.  The objective is to determine a design for which a prescribed relative eigenvalue bandgap (2$\frac{\lambda_{m+1} - \lambda_m}{\lambda_{m+1} + \lambda_m}$) is maximized, where $\lambda_m$ is the $m^{\mathrm{th}}$ eigenvalue of a certain system.  The resulting optimization problem is nonlinear, non-convex, and large-scale; nevertheless effective methods for computing solutions are developed in \cite{men2010bandgap,men2011design}.  (In fact, a more proper microstructure material design model should use binary conditions $x_i \in \{x_{\min},x_{\max}\}$ instead of the interval restrictions $x_i \in [x_{\min},x_{\max}]$; such binary conditions are typically relaxed in the bandgap optimization problems with almost no degradation in solution quality, see \cite{men2010bandgap}.)  The computed solution ${x}^{\ast}$ might not be fabricable due to small feature sizes (disconnected pixels with $x_i = x_{\max}$) or complicated material interfaces (such as roughness of boundaries).  In principle one can add constraints to ensure that features are not small and/or ensure smooth boundaries of surfaces of the material, but such an approach is decidedly unattractive as it leads to an exponential number of constraints -- and makes a complex model even more complex computationally.  Instead we proceed as follows.  First we observe that it is relatively easy in practice to use human judgment to modify a given design solution $x$ to a fabricable solution $y$ by switching the concentration of material at a relatively small number $\delta$ of pixels from $x_{\min}$ to $x_{\max}$ (or {\em vice versa}) to remove small features and/or rough material interfaces, and hopefully not degrade the objective function value too much in the process.  Here we presume that $\delta << n $.  If we anticipate that we will need to perform some sort of manual modification of a solution $x$, then it is beneficial to account for this {\em a priori} in the model specification.  This is the basic idea of {\em fabrication-adaptive optimization} which we now formally describe.\medskip

For a generic optimization problem of the form \eqref{eq_nominal}, let ${x}$ be a feasible and/or optimal solution.  The basic premise of our approach is that we will fabricate/implement some solution $y$ that is close to ${x}$ in some prescribed norm, say at most a distance $\delta$ from ${x}$ in the prescribed norm, and that such a nearby fabricable solution $y$ is very easy to determine/compute for any given solution $x$.  We construct the FA (for {\em fabrication-adaptive}) counterpart objective function $\tilde f(\cdot)$ of the original objective function $f(\cdot)$ in \eqref{eq_nominal}, as follows:
\begin{equation}\label{eq_fr-cp}
\begin{array}{ll}
\tilde{f}({x}) := & \max\limits_{y} \ \ f({y})   \\
  &  \mbox{s.t.}\ \ \ \ \|{y}-{x}\|\leq \delta  \\
  & \qquad \ \ {y}\in S \ ,
\end{array}
\end{equation}
where $\delta>0$ is the fabrication-adaptive (FA) parameter, and $\|\cdot\|$ is the prescribed norm.  Then $\tilde{f}({x})$ is a (conservative) upper bound on the objective function value of any (fabricable) solution ${y}$ whose distance from ${x}$ is at most $\delta$ in the prescribed norm.  We then construct the fabrication-adaptive optimization problem which is defined as:
\begin{equation}\label{eq_fr_opt}
\begin{array}{ll}
\tilde{z}^{\ast} = & \min\limits_{x} \ \ \tilde{f}(x)   \\
  &  \mbox{s.t.} \ \ \ x\in S \ .
\end{array}
\end{equation}
The FA optimization problem \eqref{eq_fr_opt} seeks to optimize the conservative FA counterpart $\tilde f(\cdot)$ of the original objective function $f(\cdot)$.  In this way the FA model seeks to produce a solution $x_{\mathrm{FA}}^*$ which is more adaptable to modification to a nearby  fabricable solution $y$ whose objective value $f(y)$ is not significantly degraded.\medskip

The functional form \eqref{eq_fr-cp} was first introduced by \cite{lewis2002} for the special case when $S$ is a finite-dimensional normed space (essentially $S=\mathbb{R}^n$ without loss of generality) and with the norm $\|\cdot\|$ replaced by a more general gauge function $g(\cdot)$ \citep{rock}, where it was called the ``robust regularization'' of $f(\cdot)$.  The term ``robust regularization'' is appropriate for the context given therein, which includes issues of uncertain data in the construction of $f(\cdot)$, uncertain implementation issues, and the like.  Indeed, the function $f(\cdot)$ in \cite{lewis2002} is considered broadly and so might be a constraint function or an objective function in an optimization problem, or perhaps simply a function of interest.  The name ``robust regularization'' is also suggestive of a relationship to robust optimization \citep{ben2009robust,bertsimas2011theory}, and it turns out that in the very special case when $S=\mathbb{R}^n$, the model \eqref{eq_fr-cp}-\eqref{eq_fr_opt} can be formatted as a particular instance of a robust optimization problem; this is shown herein in Appendix \ref{robustopt}.  However, when $S \ne \mathbb{R}^n$ (as one would typically expect), the connection between \eqref{eq_fr-cp}-\eqref{eq_fr_opt} and robust optimization breaks down; this is also shown in Appendix \ref{robustopt}.  \cite{lewis2009lipschitz} generalizes the definition of the robust regularization to the case when $S \ne \mathbb{R}^n$, and presents a variety of results regarding smoothness of $\tilde f(\cdot)$ and related mathematical properties.  In somewhat of a contrast, the focus of this paper is on the model \eqref{eq_fr-cp}-\eqref{eq_fr_opt} as a mechanism for fabrication-adaptive optimization; as such we rely on the premise articulated above that we will fabricate/implement some solution $y$ that is close to ${x}$ in some prescribed norm, say at most a distance $\delta$ from ${x}$ in the prescribed norm, and that such a nearby fabricable solution $y$ is very easy to determine/compute for any given solution $x$.  For this setting we expect $S \ne \mathbb{R}^n$, whereby there is no connection to robust optimization.  In the context of our intended modeling set-up, the premise of adapting the solution after-the-fact to a nearby fabricable solution, and the fact that typically $S \ne \mathbb{R}^n$, we prefer to use the term ``fabrication-adaptive optimization'' rather than ``robust regularization'' for the paradigm \eqref{eq_fr-cp}-\eqref{eq_fr_opt} as it is more aligned with aspirations of modeling, application, optimization, and computation.\medskip

The notion of fabrication adaptivity is also related to the modeling of implementation errors.  \cite{luo2003applications} and \cite{luo2004multivariate} consider such models in the context of signal processing and digital communication, where solutions are affected by errors due to discretization of the signal.  \cite{pinar2004robust} examines modeling of implementation errors in linear least-squares problems (of which signal processing is an application), and \cite{stinstra2008robust} considers implementation errors in generic optimization modeling through the lens of various types of modeling errors.  The implementation-error models developed in these works implicitly assume that the set of possible errors is independent of the solution point, therefore in this case, the modeling of implementation errors can be treated as an instance of robust optimization \citep{ben2009robust, bertsimas2011theory}.  This contrasts with fabrication adaptive optimization, as is discussed in Appendix \ref{robustopt}.\medskip

\subsection{Basic Non-Convexity Issues, and Practical and Computational Aspirations}\label{ci}

In a realistic application of the fabrication-adaptive model (\ref{eq_fr-cp})-(\ref{eq_fr_opt}) one would typically have $S \ne \mathbb{R}^n$.  However, at least from an academic perspective, the special case of $S = \mathbb{R}^n$ gives rise to interesting properties with respect to convexity and with respect to connection to the modern domain of {\em robust optimization}, see \cite{ben2009robust} and \cite{bertsimas2011theory}.  When $S = \mathbb{R}^n$, the fabrication-adaptive objective function $\tilde f(\cdot)$ is convex if $f(\cdot)$ is convex, see Proposition 3.1 of \cite{lewis2002}.  Again when $S = \mathbb{R}^n$, it is also straightforward to show that quasiconvexity is preserved as well: if $f(\cdot)$ is quasiconvex, then the fabrication-adaptive objective function $\tilde f(\cdot)$ is quasiconvex.\medskip

When $S \ne \mathbb{R}^n$, the following example shows that the FA counterpart optimization problem of a convex optimization problem need not be convex.

\begin{example}\label{peru} (A non-convex FA objective function when $f(\cdot)$ is convex.)  Let $S=\{x \in \mathbb{R}^2 : 0 \le x \le e\}$, the unit $2$-dimensional square, and consider the convex (linear) objective function $f(x):= 2x_1 + x_2$, and let $\|\cdot\|=\|\cdot\|_1$ and $\delta = 1/10$ for concreteness.  At $x^1:=(1-2\delta, 1)$ we have
$$\tilde f(x^1) = \max\{2y_1+y_2 : 0 \le y \le e, \|y-x^1\|_1 \le \delta\} = 3 -2\delta \ . $$
At $x^2:=(1, 1-3\delta)$ we have
$$\tilde f(x^2) = \max\{2y_1+y_2 : 0 \le y \le e, \|y-x^2\|_1 \le \delta\} = 3 - 2\delta \ . $$
However, at $x^3:= \frac{1}{2}x^1 + \frac{1}{2}x^2 = (1-\delta, 1-3\delta/2)$ we have:
$$\tilde f(x^3) = \max\{2y_1+y_2 : 0 \le y \le e, \|y-x^3\|_1 \le \delta\} = 3-3\delta/2 \ . $$
In this case we have $\tilde f(x^3) = 3-3\delta/2 > 3-2\delta = \max\{\tilde f(x^1), \tilde f(x^2) \}$, thus showing that $\tilde f(\cdot)$ is not even quasiconvex on the feasible region $S$, as the level set $L_\beta := \{ x \in S : \tilde f(x) \le \beta \}$ is not convex for $\beta = 3-2\delta$.\end{example}

In light of the fact that the FA optimization problem \eqref{eq_fr-cp}-\eqref{eq_fr_opt} can be non-convex, it makes most practical sense to consider using the FA modeling paradigm when the original function $f(\cdot)$ in \eqref{eq_nominal} is not convex.  (Otherwise we are doing the computational unpromising task of transforming a nominally convex problem into a non-convex problem.)  We therefore will take as given that $f(\cdot)$ is not required to be convex, and that we expect the FA optimization problem \eqref{eq_fr-cp}-\eqref{eq_fr_opt} to be non-convex.  In consideration of goals of algorithms, we aspire to compute solutions $\bar x$ of \eqref{eq_fr-cp}-\eqref{eq_fr_opt} that are either local optima or perhaps just have ``good'' objective function value $\tilde f(\bar x)$ where such ``goodness'' will of necessity be problem/context-dependent.  In order to design algorithms to solve the FA optimization problem \eqref{eq_fr-cp}-\eqref{eq_fr_opt}, we focus on two computational tasks that seem natural to require in order to design useful algorithms:  (i) computing the FA counterpart function value $\tilde f(x)$ for a given $x \in S$, and (ii) computing first-order function objects such as the gradient $\nabla \tilde f(x)$ or a (perhaps only local) subgradient of $\tilde f(x)$ or of the ``pieces'' of $\tilde f(\cdot)$ in the case when $\tilde f(\cdot)$ is the pointwise maximum of other functions, for a given $x \in S$.  For a given $x \in S$, notice from \eqref{eq_fr-cp} that computing $\tilde f(x)$ is itself generally intractable as it involves maximizing a convex function over a convex set.  Nevertheless, in many useful instances  with special structure on $f(\cdot)$, $S$, and/or $\| \cdot\|$, it will be computationally tractable to compute $\tilde{f}(\cdot)$ and $\nabla \tilde f(\cdot)$ (or other first-order information) efficiently.  Indeed, one of the main concerns of the rest of this paper is with special structures of real interest for which computation with the FA counterpart function $\tilde f(\cdot)$ is relatively efficient (Section \ref{sec_examples}), and with the practical use of the FA paradigm for solving problems that gave rise to this paradigm in the first place, namely bandgap optimization problems (Sections \ref{sec_FA-BOP} and \ref{pcd}).  In Section \ref{sec_examples} we examine FA optimization problems with certain structures of interest, mainly functions that are in turn piecewise-linear, linear fractional, piecewise linear fractional, as well as a canonical eigenvalue function.  We also propose an algorithm for FA optimization in the piecewise linear fractional case.  In Section \ref{sec_FA-BOP}, we review the class of design optimization problems known as band-gap problems that arise in engineering design optimization, and we show how the FA optimization paradigm can be applied to these problems.  In Section \ref{pcd} we present computational results from applying the FA optimization to various bandgap problems that arise in photonic crystal design, and we demonstrate that our proposed algorithm succeeds in producing much improved adaptive and fabricable solutions.

\subsection{Notation}  Let $e=(1, \ldots, 1)$ denote the vector of ones, whose dimension will be given in context.  Let $\|\cdot\|$ denote a norm on $\mathbb{R}^n$, and let $\|\cdot\|_*$ be the associated dual norm, namely $\|v\|_{\ast} := \max \{ v^T x : \|x\| \le 1\}$.  The ball of radius $\delta$ centered at $\bar x$ is denoted $B(\bar x, \delta):=\{x : \|x-\bar x\| \le \delta \}$.  Recall that a function $f(\cdot)$ on $S$ is convex if $f(\alpha x + (1-\alpha)y) \le \alpha f(x) + (1-\alpha) f(y)$ for any $\alpha \in [0,1]$ and all $x,y \in S$.  Similarly, $f(\cdot)$ on $S$ is quasiconvex if $f(\alpha x + (1-\alpha)y) \le \max\{f(x),f(y)\}$ for any $\alpha \in [0,1]$ and all $x,y \in S$, and $f(\cdot)$ is concave or quasiconcave if $-f(\cdot)$ is convex or quasiconvex, respectively.  Note that $f(\cdot)$ is quasiconvex if and only if the lower level sets of $f(\cdot)$ are convex sets, see \cite{avriel}.  If $f(\cdot)$ is convex on $S$, then $g\in\mathbb{R}^n$ is a subgradient of $f(\cdot)$ at $\hat{x}\in S$ if $f(x) \ge f(\hat{x})+ g^T(x-\hat{x})$ for all $x\in S$.  Similarly for a concave function on $S$, $g\in\mathbb{R}^n$ is a subgradient of $f(\cdot)$ at $\hat{x}\in S$ if $f(x) \le f(\hat{x})+ g^T(x-\hat{x})$ for all $x\in S$.  A function $f(\cdot)$ on $S$ is locally convex at $\hat x \in S$ if there exists some $\delta > 0$ for which $f(\cdot)$ is convex on $S \cap B(\hat x,\delta)$, and we say that $g$ is a local subgradient of $f(\cdot)$ at $\hat x$ if $f(x) \ge f(\hat{x})+ g^T(x-\hat{x})$ for all $x\in S\cap B(\hat x, \delta)$.  Similar remarks hold for local concavity and a local subgradient of a locally concave function.  Let $X, Y$ be any symmetric matrices.  We write ``$X\succeq 0$'' to denote that $X$ is symmetric and positive semidefinite, ``$X\succeq Y$'' to denote that $X-Y\succeq 0$, and ``$X\succ 0$'' to denote that $X$ is positive definite.  If $K \subset \mathbb{R}^n$ is a closed convex cone, then its dual cone $K^*$ is defined by $K^*:=\{s \in \mathbb{R}^n : s^Tx \ge 0 \ \mathrm{for~all ~} x \in K \}$.

\section{Fabrication-Adaptive Optimization Problems with Special Structures}
\label{sec_examples}
We study some FA optimization problems with special structures on $f(\cdot)$, $S$, and/or $\|\cdot\|$.
\subsection{Three Special Structures for $S=\mathbb{R}^n$}\label{eg1}

We show three classes of examples of special structures for instances where $S=\mathbb{R}^n$.  In the first class the objective function is the maximum of a finite number of affine functions:
\begin{equation}\label{eight}
f(x):= \max_{i=1,\ldots,m} b_i + (a^i)^T x \ ,
\end{equation}
and $S = \mathbb{R}^n$.  It is easy to derive the fabrication-adaptive objective function in this case:
\begin{equation}
\begin{array}{ll}
\tilde{f}(x)&  = \max\limits_{\|y-x\| \le \delta}\  \max\limits_{i=1,\ldots,m}\ b_i +  (a^i)^Ty \\[1ex]
& = \max\limits_{i=1,\ldots,m}\  \max\limits_{\|y-x\| \le \delta}\ b_i  + (a^i)^Ty \\[1ex]
& = \max\limits_{i=1,\ldots,m}\ \max\limits_{\|d\| \le \delta}\ b_i  + (a^i)^T(x+d) \\[1ex]
& = \max\limits_{i=1,\ldots,m} (b_i + \delta \|a^i\|_{\ast}) + (a^i)^Tx \ ,
\end{array}
\end{equation}
where $\|\cdot\|_{\ast}$ is the dual norm of $\|\cdot\|$.  Therefore the FA optimization problem can be written as:
\begin{equation}\label{hmm}
\min_{x\in \mathbb{R}^n} \tilde{f}(x) =  \min\limits_{x\in \mathbb{R}^n}\ \max\limits_{i=1,\ldots,m} (b_i + \delta \|a^i\|_{\ast}) + (a^i)^Tx \  .
\end{equation}
The functional form of $\tilde{f}(\cdot)$ is structurally identical to that of $f(\cdot)$, namely the maximum of $m$ linear functions, and both $f(\cdot)$ and $\tilde f(\cdot)$ are convex functions.  Let us presume that it is easy to compute the dual norm $\| a \|_*$ for any $a$.  Under this presumption the computation of $\tilde f(\cdot)$ will be as easy as that of $f(\cdot)$ and computing a subgradient of $\tilde f(\cdot)$ at a given value of $x$ will be as easy as that of $f(\cdot)$.  Furthermore, it is reasonable to expect that any algorithm for minimizing $f(\cdot)$ in \eqref{eight} should be easy to apply to solve \eqref{hmm} with similar types of computational guarantees.\medskip

The second class of examples are instances where $S=\mathbb{R}^n$ and $f(x) = \|Ax + b\|_2$ or $f(x)$ is a strictly convex quadratic function, and the prescribed norm on the space of variables $x$ is the Euclidean norm $\|x\|_2$.  In these cases, \cite{lewis2002} shows that the resulting fabrication-adaptive optimization problem \eqref{eq_fr_opt} can be modeled using semidefinite optimization.\medskip

The third class of examples are instances of the maximum eigenvalue function where $S=\mathbb{R}^n$.  Given symmetric matrices $A_0, A_1, \ldots, A_n$, let ${\cal A}(x):=A_0 + \sum_{i=1}^n A_i x_i$ and consider the maximum eigenvalue function $f(\cdot): \mathbb{R}^n \rightarrow \mathbb{R}$ given by:
$$ f(x) := \lambda_{\max}\left({\cal A}(x) \right) \ . $$
Note that $f(\cdot)$ is convex on $\mathbb{R}^n$, and $f(\cdot)$ generalizes the maximum of linear functions.  (Indeed, $f(\cdot)$ specializes to the maximum of linear functions when all matrices $A_0, A_1, \ldots, A_n$ are diagonal.)  The fabrication-adaptive counterpart function $\tilde f(\cdot)$ of $f(\cdot)$ is
\begin{equation}\label{eq_eig}
\tilde{f}({x}) = \max_{y \in S, \|y-x\| \le \delta}\ \lambda_{\max} \left({\cal A}(y)  \right) \ ,
\end{equation} and $\tilde f(\cdot)$ is also a convex function when $S=\mathbb{R}^n$ from Proposition 3.1 of \cite{lewis2002}.  However, there does not appear to be any efficient method for computing $\tilde f(\cdot)$ even in the case when $S=\mathbb{R}^n$ unless the norm $\|\cdot\|$ has very special structure.  When $\|\cdot\| = \|\cdot\|_1$, then the unit ball is the convex hull of the $2n$ signed unit vectors $e^1, \ldots, e^n, -e^1, \ldots, -e^n$, whereby $\tilde f(x)$ can be computed as:
$$\tilde f(x) = \max_{j=1, \ldots, n}\{\lambda_{\max}\left({\cal A}(x+\delta e^j)\right), \lambda_{\max}\left({\cal A}(x-\delta e^j) \right)\}\ , $$
and so is computable so long as the $2n$ largest eigenvalue problems are efficiently computable.  However, when $\|\cdot\| = \|\cdot\|_{\infty}$, it follows from \cite{Ben-tal02ontractable} that computing $\tilde f(\cdot)$ is NP-hard, and when $\|\cdot\| = \|\cdot\|_{2}$, computing $\tilde f(\cdot)$ is also NP-hard \citep{privatecom}.  Here we see that the original objective function $f(\cdot)$ involves the computation of the largest eigenvalue of a symmetric matrix, which is typically tractable; however the FA counterpart function $\tilde f(\cdot)$ is not tractable to compute when $\|\cdot\| = \|\cdot\|_p$ and $p=2$ or $p=\infty$.  Indeed, intuition suggests that only norms with a relatively small number of extreme points on their unit ball will be suitable for practical computation of the FA counterpart of the largest eigenvalue function.

\subsection{Piecewise linear convex objective and $S \ne \mathbb{R}^n$}\label{eg2}
Let the objective function be given by \eqref{eight}, i.e., the same as in Subsection \ref{eg1}, but now suppose that $S \ne \mathbb{R}^n$.  For convenience we assume in this subsection that $S$ is closed and bounded, i.e., compact.  Then we have:
\begin{equation}\label{csio}
\begin{array}{ll}
\tilde{f}(x) & = \max\limits_{y\in S, \|y-x\|\le \delta}\ \max\limits_{i = 1,\ldots,m}\ b_i + (a^i)^T y\\[1ex]
& =  \max\limits_{i = 1,\ldots,m}\ \max\limits_{y\in S, \|y-x\|\le \delta}  b_i + (a^i)^T y\\[1ex]
& = \max\limits_{i = 1,\ldots,m}\  \tilde{f}_i (x)\ ,
\end{array}
\end{equation}
where
\begin{equation}\label{csi}
\begin{array}{llc}
\tilde{f}_i(x) := &\max\limits_y\ & b_i + (a^i)^T y\\
& \mbox{s.t.} & \|y-x\| \le \delta \\
& & y \in S\ .
\end{array}
\end{equation}

\noindent We have the following result on the structure of $\tilde f(\cdot)$:

\begin{proposition}\label{propconcave} If $S$ is a convex set, then $\tilde{f}_i(\cdot): S \rightarrow \mathbb{R}$ is a concave function, $i=1, \ldots, m$, whereby $\tilde f(\cdot)$ is the pointwise maximum of concave functions.
\end{proposition}
\proof{Proof:}  Let us fix $i\in \{1, \ldots, m\}$, and let $x^1, x^2 \in S$ be given. Let $\alpha \in [0,1]$, and $x^3 = \alpha x^1 + (1-\alpha) x^2$.  Assuming for simplicity that the optimization problem defining $\tilde{f}_i(x)$ attains its optimum, let $y^j$ solve the optimization problem in the definition of $\tilde{f}_i(x^j)$ for $j=1,2$, whereby $\tilde{f}_i(x^j)=b_i + (a^i)^Ty^j$, $y^j \in S$, and $\|x^j-y^j\| \le \delta$ for $j=1,2$.  Therefore $y^3:=\alpha y^1 + (1-\alpha) y^2$ satisfies  $y^3 \in S$, and $\|x^3-y^3\| \le \delta$ and hence $y^3$ is feasible for the optimization problem in the definition of $\tilde{f}_i(x^3)$ in \eqref{csi}.  Therefore $\tilde{f}_i(x^3) \ge b_i + (a^i)^Ty^3 = \alpha( b_i + (a^i)^Ty^1) + (1-\alpha)(b_i + (a^i)^Ty^2) = \alpha \tilde{f}_i(x^1)+ (1-\alpha)\tilde{f}_i(x^2)$, and hence $\tilde{f}_i(\cdot)$ is concave on $S$.  \Halmos\medskip

It follows from Proposition \ref{propconcave} that $\tilde f(\cdot)$ does not have attractive convex structure.  Nevertheless, the computation of $\tilde f(x)$ for a given $x$ via \eqref{csio}-\eqref{csi} is a tractable convex optimization problem when $\|\cdot\|$ is the $L_2,\ L_1,\ \mbox{or } L_{\infty}$ norm, and when $S$ is polyhedral or is conveyed in a suitably easy conic form $S=\{x : b-Ax \in K\}$ for some convex cone $K$.  In these cases, computing $\tilde{f}(x)$ amounts to solving $m$ conic convex optimization problems.  This is not a particularly burdensome task if $m$ is not too large, and/or if $S$ is a relatively simple set such as a hypercube, simplex, or Euclidean ball, or more generally if $S$ is conveyed in conic form above with structure for which conic optimization can be done efficiently.\medskip

Using \eqref{csio} and \eqref{csi}, the FA optimization problem \eqref{eq_fr_opt} can therefore be written as:
\begin{equation}
\label{eq_PL_LBJ}
\begin{array}{ccccccccl}
P^{\mathrm{FA}}: \ \ \ \ \ & \tilde z^* := \ \ \ & \underset{x}{\min} & \ \  \tilde f(x) & \ \ \ \ \  = \ \ \ \ \ & \underset{x,t}{\min} &  \ \ t\\
                 &                     & \mathrm{s.t.}      & x \in S     &            &     \mathrm{s.t.} & \tilde{f}_i(x) & \le & t, \ i=1,\ldots,m\\
                 &                      &               & & &       & x              & \in & S \ . \\
\end{array}
\end{equation}
Furthermore, we know from Proposition \ref{propconcave} that $\tilde{f}_i(\cdot)$ is concave, $i=1, \dots, m$.\medskip

In light of the structure of the FA optimization problem \eqref{eq_PL_LBJ}, we consider computing first-order objects for each of the functions $\tilde{f}_i(x)$, $i=1, \ldots, m$.  Let us fix an index $i \in \{1, \ldots, m\}$.  We know from Proposition \ref{propconcave} that $\tilde f_i(\cdot)$ is concave on $S$ and hence has a subgradient for all $x\in S$.  Furthermore, there exists a set $B_i \subset S$ of measure zero such that $\tilde f_i(\cdot)$ will be differentiable for all $x \in S \setminus B_i$ (\cite{rock}, Theorem 25.5).  To see how to compute such a subgradient we appeal to duality theory and we assume that $S$ is conveyed in conic form, namely $S = \{x \in \mathbb{R}^n : b-Ax \in K \}$ where $K \subset \mathbb{R}^k$ is a closed convex cone.  For $i=1, \ldots, m$, we can re-write \eqref{csi} as the following problem $P_{ i}(x)$ :
\begin{equation}\label{csip}
\begin{array}{llc}
P_{ i}(x): \ \ \ \ \ \tilde{f}_{ i}(x) = &\max\limits_y\ & b_{ i} + (a^{ i})^T y\\
& \mbox{s.t.} & \|y-x\| \le \delta \\
& & b - Ay \in K \ ,
\end{array}
\end{equation} which can be put in conic form by defining $C:= \{(w,\alpha) : \|w\| \le \alpha\}\times K $ and re-writing $P_{ i}(x)$ as:
\begin{equation}\label{candy}
\begin{array}{llc}
P_{ i}(x): \ \ \ \ \ \tilde{f}_{ i}(x) = &\max\limits_y\ & b_{ i} + (a^{ i})^T y\\
& \mbox{s.t.} &  \left(\begin{array}{c}x \\ \delta \\ b  \end{array}\right) -\left(\begin{array}{c}I \\ 0 \\ A \end{array}\right)y \in C . \\
\end{array}\end{equation}The conic dual $D_{ i}(x)$ of $P_{ i}(x)$ can then be written as:

\begin{equation}\label{csid}
\begin{array}{llc}
D_{ i}(x): &\min\limits_\pi\ &\ b_{ i} + \delta \|a^{ i}-A^T\pi\|_* + b^T\pi + (a^{ i}-A^T\pi)^Tx \\
& \mbox{s.t.} & \pi \in K^* \ .
\end{array}
\end{equation}

\noindent We say that $S$ has a {\em Slater point} if there exists $x^0 \in S$ for which $b-Ax^0 \in \mathrm{int} K$.  The following result describes a way to compute a subgradient of $\tilde f_i(x)$:

\begin{proposition}\label{supergradient} Let $i\in \{1, \ldots, m\}$ be given.  Suppose that $S$ has a Slater point, and suppose $x \in S$.  Then $D_{ i}(x)$ attains its optimum at some $\pi_i^*$ with no duality gap, and furthermore \begin{equation}\label{super} p_i:=p_i(x):=a^{ i} - A^T\pi_i^* \end{equation} is a subgradient of $\tilde f_i(\cdot)$ at $x$.  Furthermore, there is a set $B_i \subset S$ of measure zero for which it holds that $p_i(x)$ is uniquely defined and hence $\nabla \tilde f_i(x) = p_i(x)$ for all $x \in S \setminus B_i$.
\end{proposition}

\proof{\bf Proof:}  Let us fix $i \in \{1, \ldots, m\}$ and consider the duality paired problems $P_{ i}(x)$ and $D_{ i}(x)$. It follows from standard duality theory that there will be no duality gap and the dual problem will attain its optimum under the condition that the primal has a Slater point, namely a point $y$ for which $\|y-x\| < \delta$ and $b-Ay \in \mathrm{int} K$; see \cite{duffin} or \cite{borweinlewis} for a more modern treatment of conic duality.  Let $x^0$ be a Slater point of $S$, whereby $b-Ax^0 \in \mathrm{int} K$.  Since $x \in S$ by supposition, it follows that $y(\varepsilon):= \varepsilon x^0 + (1-\varepsilon) x$ is a Slater point of the feasible region of \eqref{candy} for all $\varepsilon >0$ and sufficiently small.  It then follows that $D_i(x)$ attains its optimum with no duality gap, and it follows from the formulation of $D_i(x)$ that $\tilde f_i(x)$ is the pointwise minimum of affine functions, whose linear terms are of the form $(a^{ i}-A^T\pi)$ for $\pi \in K^*$.  It then follows directly from convexity arguments that $a^{ i} - A^T\pi_i^* $ is a subgradient of $\tilde f_i(\cdot)$ at $x$. Furthermore, it follows from \cite{rock} (Theorem 25.5) that there exists a set $B_i \subset S$ of measure zero such that $\tilde f_i(\cdot)$ is differentiable for all $x \in S \setminus B_i$, and hence $\nabla \tilde f_i(x) = p_i(x)$ for all $x \in S \setminus B_i$.  \Halmos\medskip

The computational viability of solving the dual problem $D_{ i}(x)$ must of necessity presume that the dual norm $\|\cdot\|_*$ can be suitably treated in the objective function of $D_{ i}(x)$.  Of course, when the norm $\|\cdot\|$ can be described with linear inequalities or second-order cone constraints, then solving $D_{ i}(x)$ is all the more easy.  For example, when $\|\cdot\|$ is the $L_1$- or $L_\infty$-norm, then $D_{ i}(x)$ can be easily represented as a linear programming problem.  When $\|\cdot\|$ is the $L_2$-norm or other quadratic norm of the form $\sqrt{x^TQx}$ for $Q \succ 0$, then $D_{ i}(x)$ can be represented as a second-order cone problem using a standard transformation, see \cite{vb}.\medskip

In Section \ref{eg4} we will present an algorithm for solving \eqref{eq_PL_LBJ} that is based on the scheme of sequentially solving the first-order approximation of \eqref{eq_PL_LBJ} at a given point $\hat x \in S$, namely:
\begin{equation}
\label{eq_PL_GF}
\begin{array}{llrcl}
P^{\hat x}: & \underset{x,t}{\min} & t\\
& {s.t.} &\ \tilde{f}_i(\hat x) + \nabla \tilde f_i(\hat x)^T(x - \hat x) &\le& t \ , \qquad i=1,\ldots,m\\
&& x &\in& S \ ,
\end{array}
\end{equation} where the values $\tilde {f}_i(\hat x)$ and $\nabla \tilde{f}_i(\hat x)$, $i=1, \ldots, m$, are computed via Propositions \ref{propconcave} and \ref{supergradient}, respectively.\medskip

\subsection{Linear fractional objective and $S \ne \mathbb{R}^n$}\label{eg3}
Let us now consider the case when the objective function is linear fractional:
\begin{equation} f(x) := \frac{a^Tx + g}{c^Tx + h} \ , \end{equation} and suppose that $S \ne \mathbb{R}^n$ and we impose the condition that $S$ is compact and convex, and $c^Tx + h > 0$ for all $x \in S$.  In this case $f(\cdot)$ is quasilinear on $S$, i.e., it is both quasiconvex and quasiconcave on $S$.  We can write the fabrication-adaptive objective function as:

\begin{equation}\label{91}
\begin{array}{llc}
\tilde f(x) := &\max\limits_y\ & \displaystyle\frac{a^Ty + g}{c^Ty + h}\\[2ex]
& \mbox{s.t.} & \|y-x\| \le \delta \\[1ex]
& & y \in S\ .
\end{array}
\end{equation}

\noindent We have the following result on the structure of $\tilde f(\cdot)$:

\begin{proposition}\label{propconcave2} $\tilde f(\cdot): S \rightarrow \mathbb{R}$ is a quasiconcave function on $S$.
\end{proposition}
\proof{\bf Proof:}  Let $x^1, x^2 \in S$ be given, let $\alpha \in [0,1]$, and let $x^3 = \alpha x^1 + (1-\alpha) x^2$.  Let $y^j$ solve the optimization problem in the definition of $\tilde f(x^j)$ in \eqref{91} for $j=1,2$, whereby $\|x^j-y^j\| \le \delta$, $y^j \in S$, and  $\frac{a^Ty^j + g}{c^Ty^j + h}=\tilde f(x^j)$ for $j=1,2$.  Define $\beta := \min\{\tilde f(x^1), \tilde f(x^2)\}$.  Then the equations $\tilde f(x^j)=\frac{a^Ty^j + g}{c^Ty^j + h}$ for $j=1,2$ and the definition of $\beta$ implies that $a^Ty^j + g = \tilde f(x^j) [c^Ty^j + h] \ge  \beta [c^Ty^j + h]$ for  $j=1,2$.  Also $y^3:=\alpha y^1 + (1-\alpha) y^2$ satisfies  $y^3 \in S$ and $\|x^3-y^3\| \le \delta$ and hence is feasible for the optimization problem in the definition of $\tilde f(x^3)$ in \eqref{91}.  Furthermore $a^Ty^3 + g \ge \beta [c^Ty^3 + h]$, and hence $\tilde f(x^3) \ge \beta = \min\{\tilde f(x^1), \tilde f(x^2)\}$, whereby $\tilde f(\cdot)$ is quasiconcave on $S$.  \Halmos\medskip

The next result shows that $\tilde f(x)$ is the optimal objective function value of a convex optimization problem involving $S$ and a norm constraint.  In what follows we use the notation that $\theta S$ is the scaling of $S$ by the constant $\theta$ for $\theta \ge 0$.  (When $\theta =0$, it is customary to define $\theta S$ to be the recession cone of $S$; however since we assume here that $S$ is bounded we have for $\theta = 0$ that $\theta S := \{0\}$ under either definition.)

\begin{proposition}\label{92} If $x\in S$, then $\tilde f(x)$ is the optimal objective value of the following convex optimization problem:
\begin{equation}\label{93}
\begin{array}{llrll}
\tilde f(x) := &\max\limits_{\bar y, \theta}\ & a^T{\bar y} + g\theta \\[1ex]
& \mbox{s.t.} & \|\bar y-\theta x\| &\le& \theta\delta \\[1ex]
& & c^T{\bar y} + h\theta & = &1 \\[1ex]
& & \bar y & \in & \theta S\\[1ex]
& & \theta &\ge &0 \ .
\end{array}
\end{equation}  If $(\bar y, \theta)$ is feasible and/or optimal for \eqref{93}, then $y = \bar y/\theta$ is feasible and/or optimal for \eqref{91}, respectively.
\end{proposition}

\proof{\bf Proof:}  Problem \eqref{93} is just a standard transformation of the linear fractional optimization problem \eqref{91} via homogenization (see \cite{charnes1962plf} or \cite{craven1973dfl}).  Note that we cannot have $\theta=0$ in \eqref{93}, as this would imply that $\bar y = 0$ via the norm constraint, which would imply that $c^T\bar y + h \theta =0 \ne 1$, which is a contradiction.  Thus division by zero in the transformation $y=\bar y/\theta$ cannot occur.  \Halmos\medskip

The computational viability of solving \eqref{93} will depend on (among other things) the ability to conveniently work with the scaling $\theta S$.  When $S$ is conveyed in conic linear form as in Section \ref{eg2}, it follows for $\theta \ge 0$ that $\theta S = \{ x \in \mathbb{R}^n : \theta b - Ax \in K \}$ , in which case the scaling results in no loss of generality of the representation of the feasible region.\medskip

Let us now turn to the computation of first-order objects related to $\tilde f(\cdot)$.  Let us further assume that $S$ is conveyed as a system of linear inequalities, so that $S = \{ x \in \mathbb{R}^n : b - Ax \ge 0 \}$, i.e., $K=\mathbb{R}^k_+$ in the conic linear representation.  We then can write \eqref{93} as:

\begin{equation}\label{94}
\begin{array}{llrll}
\tilde f(x) := &\max\limits_{\bar y, \theta}\ & a^T{\bar y} + g\theta \\[1ex]
& \mbox{s.t.} & \|\bar y-\theta x\| &\le& \theta\delta \\[1ex]
& & c^T{\bar y} + h\theta  &=& 1 \\[1ex]
& & \theta b - A \bar y &\ge& 0\\[1ex]
& & \theta &\ge& 0 \ .
\end{array}
\end{equation}
Let us further restrict our attention to the case when the norm $\|\cdot\|$ is representable with linear inequalities, as in the $L_1$- or $L_\infty$-norm.  For concreteness let us examine the case when $\|\cdot\|$ is a weighted $L_1$-norm with weights $w:=(w_1, \ldots, w_n) >0$, i.e., $\|x\|:=\sum_{j=1}^n w_j|x_j|$.  Then \eqref{94} can be represented as the following linear optimization problem:

\begin{equation}\label{95}
\begin{array}{llrllr}
\tilde f(x) := &\max\limits_{\bar y, \theta, q}\ &\quad a^T{\bar y} + g\theta \\[1ex]
& \mbox{s.t.} & \bar y - \theta x &\le& q &  \ \ \ \ \ (\pi_1)\\[1ex]
& & -\bar y + \theta x &\le& q & (\pi_2)\\[1ex]
& & w^Tq &\le& \theta \delta & (\gamma)\\[1ex]
& & c^T{\bar y} + h\theta  &=& 1 & (\tau)\\[1ex]
& & \theta b - A \bar y &\ge& 0 & (\lambda)\\[1ex]
& & \theta &\ge 0& \ , & (\kappa)
\end{array}
\end{equation} where for future reference we assign names of linear optimization dual variables for each of the constraint systems above.  We have the following result concerning the computation of the gradient of $\tilde f(\cdot)$:

\begin{proposition}\label{supergradient2} There is a set $B \subset S$ of measure zero which makes the following hold:  suppose $x \in S\setminus B$, and let $(\bar y^*,\theta^*, q^*)$ solve \eqref{95} and let $(\pi_1^*, \pi_2^*, \gamma^*, \tau^*, \lambda^*, \kappa^*)$ be optimal dual variables, and define:
\begin{equation}\label{super2} p:=p(x):= \theta^*(\pi^*_1-\pi^*_2) \ .  \end{equation} Then $\nabla \tilde f(x) = p(x)$.
\end{proposition}

Below we present a proof of Proposition \ref{supergradient2}.  This proof relies on \eqref{95} being a linear optimization problem.  This will be the case whenever $S$ is conveyed via linear inequalities, and whenever the norm $\|\cdot\|$ is polyhedral, or to be more exact, whenever the norm level set constraint ``$\|v\| \le t$'' in variables $v,t$ can be conveniently represented via linear inequalities.  While we specifically worked with the weighted $L_1$-norm in \eqref{95}, there is no loss of generality in working with a weighted $L_\infty$-norm or other polyhedral norm.\medskip

As \eqref{95} is a linear program parameterized by $x$, it would be convenient to prove Proposition \ref{supergradient2} by invoking a standard right-hand-side sensitivity analysis result on parametric linear programming. However, notice that the parameter $x$ appears in the left-hand-side of the first two constraints of \eqref{95}, and this dependence of the constraint-matrix coefficients of \eqref{95} on the parameter $x$ appears to be structural, i.e., we see no way to remove it by a simple change of variable.  Thus to prove Proposition \ref{supergradient2} we will invoke the following result concerning changes in data coefficients in linear programming:

\begin{theorem}\label{sa}  Consider the following primal and dual pair of linear optimization problems:
\begin{equation}
\begin{array}{lrlclll}
P({x}): & f(x)=\underset{y}{\min} & c_{x}^T y &\qquad& D(x): & \underset{z}{\max} & b_{x}^T z\\[1ex]
&s.t.&A_{x} y = b_{x}, & & &s.t.&A^T_{x} z \le c_{x},\\[1ex]
& & y \ge 0.& & & & z \mbox{ free}.
\end{array}
\end{equation}
Suppose $A_{x} = A_0 + \sum_{k=1}^{n_x} A_k x_k$, $b_{x} = b_0 + \sum_{k=1}^{n_x} b_k x_k$, and $c_{x} = c_0 + \sum_{k=1}^{n_x} c_k x_k$, where $x = (x_1, x_2,\ldots,x_{n_x})$ are parameters that determine the data $A_{x}, b_{x}, c_{x}$ of the linear program $P(x)$.  Let $S$ denote the subset of $\mathbb{R}^{n_x}$ for which $P(x)$ has an optimal solution.  Then there exists a set $B \subset S$ of measure zero which makes the following hold: if $\hat{x} \in S \setminus B$, and $y^{\ast}, z^{\ast}$ are optimal solutions to $P(\hat{x})$ and $D(\hat{x})$, then
\begin{equation}
\left.\frac{\partial f(x)}{\partial x_k} \right\vert_{x_k = \hat{x}_k} = c_k^T y^{\ast} + b_k^T z^{\ast} - (z^{\ast})^T A_k y^{\ast}.
\end{equation}
\end{theorem}
Theorem \ref{sa} follows as the multivariate extension of the case when $n_x=1$ using rational functions (see Lemma 1 and Theorem 2 of \cite{freund1985postoptimal}).\Halmos \medskip

\proof{\bf Proof of Proposition \ref{supergradient2}:} The proof follows by applying Theorem \ref{sa} in the linear optimization problem \eqref{95}.  Considering the $k^{\mathrm{th}}$ component of $x$ we have from Theorem \ref{sa} that
$$ \left.\frac{\partial \tilde f(x)}{\partial x_k} \right\vert_{x_k} = \theta^*(\pi^*_1)_k - \theta^*(\pi^*_2)_k \ , $$ except possibly on a set $B$ of measure zero, which proves the result.\Halmos\medskip

\subsection{Piecewise-linear fractional objective and $S \ne \mathbb{R}^n$}\label{eg4}

Let us now consider the case when the objective function is piecewise-linear fractional:
\begin{equation}\label{122} f(x) := \max_{i=1, \ldots, m}\frac{(a^i)^Tx + g_i}{(c^i)^Tx + h_i} \ , \end{equation} and we impose the conditions that $S$ is compact and convex, and $(c^i)^Tx + h_i > 0$ for all $x \in S$ and for all $i=1, \ldots, m$.

Similar to Section \ref{eg3}, it holds that $f(\cdot)$ is quasiconvex on $S$ (but not quasiconcave).  We can write the FA objective function as:

\begin{equation}\label{csi2}
\begin{array}{ll}
\tilde{f}(x) & = \max\limits_{y\in S, \|y-x\|\le \delta}\ \max\limits_{i = 1,\ldots,m} \displaystyle \frac{(a^i)^Ty + g_i}{(c^i)^Ty + h_i}\\[2ex]
& =  \max\limits_{i = 1,\ldots,m}\ \max\limits_{y\in S, \|y-x\|\le \delta}  \displaystyle\frac{(a^i)^Ty + g_i}{(c^i)^Ty + h_i}\\[1ex]
& = \max\limits_{i = 1,\ldots,m}\ \tilde{f}_i (x)\ ,
\end{array}
\end{equation}
where
\begin{equation}\label{csi3}
\begin{array}{llc}
\tilde{f}_i(x) := &\max\limits_y\ & \displaystyle\frac{(a^i)^Ty + g_i}{(c^i)^Ty + h_i}\\[2ex]
& \mbox{s.t.} & \|y-x\| \le \delta \\
& & y \in S\ .
\end{array}
\end{equation}

\noindent The following result on the structure of $\tilde f(\cdot)$ is evident from Proposition \ref{propconcave2}:

\begin{proposition}\label{propconcave3} $\tilde{f}_i(\cdot): S \rightarrow \mathbb{R}$ is a quasiconcave function, $i=1, \ldots, m$, whereby $\tilde f(\cdot)$ is the pointwise maximum of quasiconcave functions.
\end{proposition}
\proof{\bf Proof: } The proof is an immediate consequence of Proposition \ref{propconcave2}.  \Halmos\medskip

Paralleling results in Section \ref{eg3}, we have the following result on the computation of $\tilde f(\cdot)$, which shows that $\tilde f(\cdot)$ is computable by solving $m$ convex optimization problems.

\begin{proposition}\label{97} If $x\in S$, then $\tilde f(x) = \max\limits_{i=1, \ldots, m} \tilde{f}_i(x)$ where for $i=1, \ldots, m$, $\tilde{f}_i(x)$ is the optimal objective function value of the following convex optimization problem:
\begin{equation}\label{98}
\begin{array}{llc}
P_i(x): \ \ \ \ \ \ \tilde{f}_i(x) := &\max\limits_{\bar y, \theta}\ & (a^i)^T{\bar y} + g_i\theta \\[1ex]
& \mbox{s.t.} & \|\bar y-\theta x\| \le \theta\delta \\[1ex]
& & (c^i)^T{\bar y} + h_i\theta  = 1 \\[1ex]
& & \bar y \in \theta S\\[1ex]
& & \theta \ge 0 \ .
\end{array}
\end{equation}  For each $i=1, \ldots, m$, if $(\bar y, \theta)$ is feasible and/or optimal for $ P_i(x)$ in \eqref{98}, then $y = \bar y/\theta$ is feasible and/or optimal for \eqref{csi3}, respectively.
\end{proposition}

\proof{\bf Proof: }  Proposition \ref{97} is essentially a restatement of Proposition \ref{92} for each of the $i=1, \ldots, m$ pieces $\tilde f_i(\cdot)$ of $\tilde f(\cdot)$. \Halmos \medskip

As in Section \ref{eg2}, using \eqref{csi2} and \eqref{csi3} the FA optimization problem \eqref{eq_fr_opt} can be written as:
\begin{equation}
\label{eq_PL_FR0}
\begin{array}{ccccccccl}
P^{\mathrm{FA}}: \ \ \ \ \ & \tilde z^* := \ \ \ & \underset{x}{\min} & \ \  \tilde f(x) & \ \ \ \ \  = \ \ \ \ \ & \underset{x,t}{\min} &  \ \ t\\
                 &                     & \mathrm{s.t.}      & x \in S     &            &     \mathrm{s.t.} & \tilde{f}_i(x) & \le & t, \ i=1,\ldots,m\\
                 &                      &               & & &       & x              & \in & S \ . \\
\end{array}
\end{equation}
Furthermore, we know from Proposition \ref{propconcave3} that $\tilde{f}_i(\cdot)$ is quasiconcave, $i=1, \dots, m$.\medskip

Let $\hat x \in S$ be a given point.  In light of the structure of the fabrication-adaptive optimization problem \eqref{eq_PL_FR0}, we consider computing first-order objects for each of the functions $\tilde{f}_i(x)$, $i=1, \ldots, m$.  Similar to Section \ref{eg3}, we assume that $S$ is conveyed as a system of $k$ linear inequalities, namely $S = \{ x \in \mathbb{R}^n : b - Ax \ge 0 \}$, i.e., $K=\mathbb{R}^k_+$ in the conic linear representation, and in particular we examine the case when the prescribed norm is the weighted $L_1$-norm with weights $w:=(w_1, \ldots, w_n) >0$.  Then for $i=1, \ldots, m$, problem \eqref{98} can be represented as the following linear optimization problem:

\begin{equation}\label{eq_PLi}
\begin{array}{llrlr}
P_i(x): \ \ \ \ \ \ \tilde{f}_i(x) := &\max\limits_{\bar y_i, \theta_i, q_i}\ &\quad  (a^i)^T{\bar y_i} + g_i\theta_i \\[1ex]
& \mbox{s.t.} & \bar y_i - \theta_i x & \le q_i & \ \ \ \ \  \ \ \ \ \ (\pi_{i,1})\\[1ex]
& & -\bar y_i + \theta_i x & \le q_i & (\pi_{i,2})\\[1ex]
& & w^Tq_i & \le \theta_i \delta & (\gamma_i)\\[1ex]
& & (c^i)^T{\bar y_i} + h_i\theta_i &  = 1 & (\tau_i)\\[1ex]
& & \theta b - A \bar y_i & \ge 0 & (\lambda_i)\\[1ex]
& & \theta_i & \ge 0 \ , & (\kappa_i)
\end{array}
\end{equation} where for future reference we assign names of linear optimization dual variables for each of the constraint systems above.  The analogous result of Proposition \ref{supergradient2} is:

\begin{proposition}\label{supergradient3} Let $i\in \{1, \ldots, m\}$ be given.  There is a set $B_i \subset S$ of measure zero which makes the following hold: suppose $x \in S \setminus B_i$.  Let $(\bar y_i^*,\theta_i^*, q_i^*)$ solve $P_{i}(x)$ of \eqref{eq_PLi} and let $(\pi_{i,1}^*, \pi_{i,2}^*, \gamma_i^*, \tau_i^*, \lambda_i^*, \kappa_i^*)$ be optimal dual variables, and define:
\begin{equation}\label{eq_PLi_grad} p_i:=p_i(x):= \theta_i^*(\pi^*_{i,1}-\pi^*_{i,2}) \ .  \end{equation}
Then $\nabla \tilde f_i(x) = p_i(x)$.
\end{proposition}

\proof{\bf Proof:} The result follows directly from Proposition \ref{supergradient2}.  \Halmos\medskip

Based on Propositions \ref{97} and \ref{supergradient3}, we propose an algorithm for solving the FA optimization problem $P^{\mathrm{FA}}$ in \eqref{eq_PL_FR0}.  Let $\hat x \in S$ be a given point.  We sequentially solve the first-order approximation of \eqref{eq_PL_FR0} based on the point $\hat x$, namely:
\begin{equation}
\label{eq_PL_FR1}
\begin{array}{llrcl}
P^{\hat x}: & \underset{x,t}{\min} & t\\
& {s.t.} & \tilde{f}_i(\hat x) + \nabla \tilde f_i(\hat x)^T(x - \hat x) &\le& t \ , \qquad i=1,\ldots,m\\
&& x &\in& S \ ,
\end{array}
\end{equation} where the values $\tilde {f}_i(\hat x)$ and $\nabla \tilde{f}_i(\hat x)$, $i=1, \ldots, m$ are computed via Propositions \ref{97} and \ref{supergradient3}, respectively.  This leads to the sequential linear optimization scheme described in Table \ref{tb_PL}, which we refer to as Algorithm FA (for Fabrication-Adaptivity).\medskip

\begin{table}
\caption{Algorithm for fabrication-adaptive optimization problem when $f(\cdot)$ is piecewise linear fractional and $S$ is given by linear inequalities.\label{tb_PL}}
\begin{center}
\begin{tabular}{|ll|}
\hline
& \textbf{Algorithm FA when $f(\cdot)$ is Piecewise Linear Fractional Problem}\\\hline
{\bf Step 1.} & Start with initial guess $\hat x := x^0$ and tolerance $\epsilon_{\mbox{tol}}$\\  & \\
{\bf Step 2.} & For each $i =1, \dots, m$, do:\\
              & \quad Compute function value $\tilde{f}_i(\hat{x})$ via \eqref{eq_PLi}\\
              & \quad Compute first-order information $p_i(\hat{x})$ via \eqref{eq_PLi_grad}\\ & \\
{\bf Step 3.} & Form the linear optimization problem $P^{\hat x}$ in \eqref{eq_PL_FR1} \\ & \\
{\bf Step 4.} & Solve $P^{\hat x}$ for an optimal solution $x^*$\\ & \\
{\bf Step 5.} & If $\|x^* - \hat x\| \le \epsilon_{\mbox{tol}}$, stop.\\
              & Else update $\hat x \leftarrow x^*$ and go to {\bf Step 2.}\\
\hline
\end{tabular}
\end{center}
\end{table}

Note that Algorithm FA is designed for the case when $f(\cdot)$ is piecewise linear fractional \eqref{122}.  When the denominators in the linear fractional forms are all equal to $1$, i.e., $c^i=0$ and $h_i=1$ for $i=1, \ldots, m$, then it follows from Proposition \ref{propconcave} that $\tilde{f}_i(\cdot)$ is concave on $S$, whereby for $\hat x \in S$ we have:
$$\tilde{f}_i(x) \le \tilde{f}_i(\hat x) + \nabla \tilde f_i(\hat x)^T(x - \hat x) \ \mathrm{for~all~} s\in S, \ i=1, \ldots, m \ .  $$  This in turn implies that if $(x,t)$ is feasible for $P^{\hat x}$ \eqref{eq_PL_FR1}, then $(x,t)$ is also feasible for $P^{\mathrm{FA}}$ \eqref{eq_PL_FR0}, whereby the optimal value of $P^{\hat x}$ will always be an upper bound on the value of $P^{\mathrm{FA}}$ in this case.\medskip

When we discuss bandgap optimization problems in Section \ref{sec_FA-BOP} (for which the concept of fabrication adaptivity was originally inspired), we will show that the original bandgap optimization problem can be cast as an instance of the piecewise linear fractional optimization problem \eqref{122}, and that Algorithm FA can therefore be used to solve the FA optimization problem associated with this problem.  We will show computational results for a particularly useful bandgap optimization problem, namely the photonic crystal design problem, in Section \ref{sec_res}.

\subsection{A Very Special Piecewise Linear Fractional Problem}\label{eg5}
Let us now consider the following very special piecewise linear fractional objective function:
\begin{equation}\label{sandy}
f(x) := \frac{\underset{i \in \cal I}{\max} ((a^i)^T x + g_i) - \underset{j \in \cal J}{\min} ((c^j)^T x + h_j)}{\underset{i \in \cal I}{\max} ((a^i)^T x + g_i) + \underset{j \in \cal J}{\min} ((c^j)^T x + h_j)} \ .
\end{equation}
Similar in spirit to Section \ref{eg4}, we suppose that $S$ is compact and convex, and we impose the condition that $(a^i)^T x + g_i > 0$ and $(c^j)^T x + h_j>0$ for all $x \in S$, $i\in \cal I$, and $j \in \cal J$.  Optimization problems with this structure arise naturally and often in bandgap optimization applications, which will be discussed in Section \ref{sec_FA-BOP}.  (Indeed, bandgap optimization problems gave rise to our interest in this particular structure to begin with.)

In order to analyze the properties of $f(\cdot)$ as well as the fabrication-adaptive objective $\tilde f(\cdot)$ we will use the following result.

\begin{proposition}\label{separation} Suppose that $(a^i)^T x + g_i > 0$ and $(c^j)^T x + h_j>0$ for all $x \in S$, $i\in \cal I$, and $j \in \cal J$.  Then
\begin{equation}\label{121} f(x) := \frac{\underset{i \in \cal I}{\max} ((a^i)^T x + g_i) - \underset{j \in \cal J}{\min} ((c^j)^T x + h_j)}{\underset{i \in \cal I}{\max} ((a^i)^T x + g_i) + \underset{j \in \cal J}{\min} ((c^j)^T x + h_j)} = \underset{{i \in \cal I}, {j \in \cal J}}{\max}\frac{((a^i)^T x + g_i) -((c^j)^T x + h_j) }{((a^i)^T x + g_i) + ((c^j)^T x + h_j)} \ . \end{equation}
\end{proposition}

\proof {\bf Proof:} To ease the notational burden let $U_i = (a^i)^T x + g_i$ for $i\in \cal I$ and $L_j=(c^j)^T x + h_j$ for $j \in \cal J$.  Let us also define the general function $\phi(U,L):=\frac{U-L}{U+L}$, and notice that $\phi(\cdot,\cdot)$ is increasing in $U$ and decreasing in $L$ for $U>0$ and $L>0$.  For a given value of $x\in S$, let us consider the left side and the right side of the equality \eqref{121} we need to prove.  Clearly the right side is at least as large as the left side.  Now suppose that $(\widehat i,\widehat j)$ is a pair of indices that attains the maximum in the right side of \eqref{121}.  Then we have from the monotonicity of $\phi(U,L)$ that
$$\frac{{\max\limits_{i \in \cal I}}\ U_i - {\min\limits_{j \in \cal J}}\ L_j}{{\max\limits_{i \in \cal I}}\ U_i + {\min\limits_{j \in \cal J}}\ L_j} \ge \frac{U_{\widehat i} - {\min_{j \in \cal J}}\ L_j}{U_{\widehat i} + {\min_{j \in \cal J}}\ L_j} \ge
\frac{U_{\widehat i} - L_{\widehat j}}{U_{\widehat i} + L_{\widehat j}} \ ,
$$ thus showing that the left side of \eqref{121} is at least as large as the right side, completing the proof. \Halmos\medskip

Proposition \ref{separation} shows that $f(\cdot)$ can alternatively be rewritten as the maximum of $m:=|{\cal I}| \cdot |{\cal J}|$ linear fractional functions (where $|\cal K|$ denotes the cardinality of the set $\cal K$), and so is an instance of the format \eqref{122}, and hence all of the results of Section \ref{eg4} apply herein.\medskip

Indeed, we can write the FA objective function as:
\begin{equation}\label{desk}
\begin{array}{ll}
\tilde{f}(x) & = \max\limits_{y\in S, \|y-x\|\le \delta}\ \displaystyle\frac{\max\limits_{{i \in \cal I}} ((a^i)^T y + g_i) - \min\limits_{j \in \cal J} ((c^j)^T y + h_j)}{\max\limits_{i \in \cal I} ((a^i)^T y + g_i) + \min\limits_{j \in \cal J} ((c^j)^T y + h_j)}  \\[1ex]
& = \max\limits_{y\in S, \|y-x\|\le \delta}\ \max\limits_{{i \in \cal I},{j \in \cal J}} \displaystyle\frac{((a^i)^T y + g_i) - ((c^j)^T y + h_j)}{((a^i)^T y + g_i) + ((c^j)^T y + h_j)} \\[1ex]
& = \max\limits_{{i \in \cal I},{j \in \cal J}}\ \max\limits_{y\in S, \|y-x\|\le \delta} \displaystyle\frac{(a^i - c^j)^T y + (g_i - h_j)}{(a^i + c^j)^T y + (g_i + h_j)} \\[2ex]
& =  \max\limits_{{i \in \cal I},{j \in \cal J}}\  \tilde{f}_{i,j} (x)\ ,
\end{array}
\end{equation}
where for $(i,j) \in {\cal I} \times {\cal J}$ we have:
\begin{equation}\label{fijast}
\begin{array}{llc}
\tilde{f}_{i,j}(x) := &\max\limits_y\ &  \displaystyle\frac{(a^i - c^j)^T y + (g_i - h_j)}{(a^i + c^j)^T y + (g_i + h_j)} \\[2ex]
& \mbox{s.t.} & \|y-x\| \le \delta \\[1ex]
& & y \in S\ ,
\end{array}
\end{equation}
which has the exact same format as \eqref{csi2} and \eqref{csi3} with a finite index set given by all pairs $(i,j) \in {\cal I} \times {\cal J}$.  Therefore $\tilde{f}_{i,j}(\cdot)$ is quasiconcave and $\tilde f(\cdot)$ is the pointwise maximum of quasiconcave functions (Proposition \ref{propconcave3}), the computation of $\tilde{f}_{i,j}(x)$ for a given $x \in S$ is the solution of the convex optimization problem described in Proposition \ref{97}, and the computation of $\nabla \tilde{f}_{i,j}(x)$ for a given $x$ is obtained from dual variables as described in Proposition \ref{supergradient3}.  Furthermore, the algorithm presented in Section \ref{eg4} can be applied to solve the FA optimization problem derived from the original objective function given in \eqref{sandy}.

\section{Fabrication Adaptivity for Bandgap Optimization Problems}
\label{sec_FA-BOP}
The motivation for developing the fabrication adaptivity paradigm stemmed from work on bandgap design optimization, and more specifically on photonic crystal design optimization. The works of \cite{cox2000bso,kao2005mbg}, and \cite{men2010bandgap,men2011design} contain methods to optimize bandgaps for this class of problems, but of necessity none of these works address issues of fabricability.  The chief goal of this paper is to construct and solve fabrication-adaptive optimization problems for this class of problems.  In this section we first review bandgap optimization problems in general, and we present the class of bandgap optimization models used in the previous work.  We then show how to apply the fabrication-adaptive formulation to bandgap optimization problems.  The presentation herein is at a high level, and properties are stated as summary results of previous work without proofs; for a more detailed presentation we refer the interested reader to \cite{men2010bandgap}.

\subsection{Bandgap Optimization Problems}
A {\em bandgap} is a concept that arises in many engineering applications. In semiconductor physics, an electron can sometimes transition from one energy state to another by a change in crystal momentum. An {\em energy bandgap} thus denotes a range of energy states that the electrons are forbidden to occupy despite any change in momentum (in the absence of any external excitation). Analogously, in photonic crystals (periodic optical nanostructures), photons can behave as waves, and propagate with certain frequencies through the bulk material at admissible wavevectors. A {\em frequency bandgap} is defined as the range of disallowed frequencies of the photons; if a photon is traveling according to a given wavevector, it will get attenuated within the crystal if it is propagating at any frequency within the frequency bandgap.  The energy bandgap phenomenon has been used in applications such as insulators, laser diodes, solar cells, etc., while the frequency bandgap phenomenon has been used in applications such as frequency filters, waveguides, and optical buffers.\medskip

\begin{figure}[t]
\begin{center}
\includegraphics[scale=0.8]{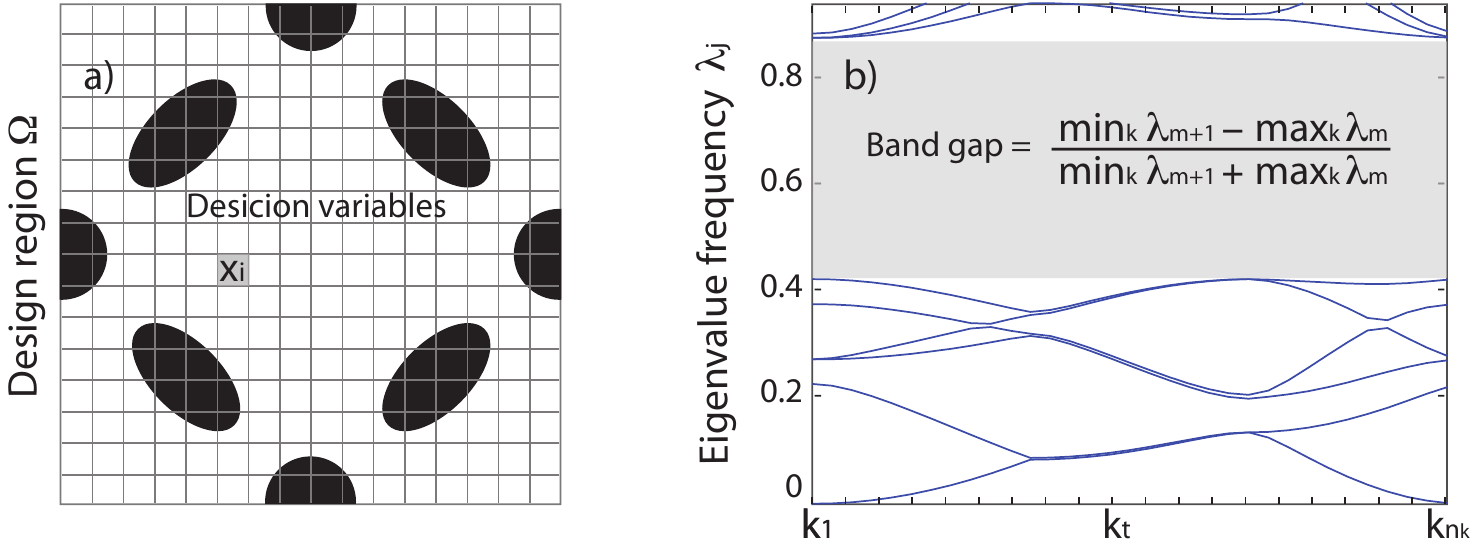}
\caption{Bandgap optimization problem. The left figure is the schematic representation of the design problem.  The design region is represented by  piecewise constant values $x_i$ of pixel (or voxel) $i$ for each of $i=1, \ldots,n_x$ pixels, where $x_i$ is the value of the material property of design interest (such as permittivity, Poisson's ratio, Young's modulus, etc.).  The right figure is the band diagram:  in this example $\cal Q$ is the interval $[k_1, k_{n_k}]$, which is discretized into a finite number of components $\{k_1, \ldots, k_{n_k}\}$.  Each curve shows an eigenvalue parameterized over the elements of $\cal Q$.  The band shown in the figure is the difference between the $6^{\mathrm{th}}$ and $7^{\mathrm{th}}$ eigenvalues, and portrays the numerator of the objective function in \eqref{eqBGopt}.  The objective is to determine the values of the design variables $x$ for which the resulting bandgap between these consecutive eigenvalues is largest.}
\label{fig_BGscheme}
\end{center}
\end{figure}

Bandgap optimization is the process of designing the composition and structure of a material to maximize a specific bandgap.  The bandgap optimization problem is generally written as the following nonlinear non-convex eigenvalue-constrained optimization problem:
\begin{equation}\label{eqBGopt}
\begin{array}{lcll}
P: &  \underset{ x\in S}{\max} & \ \ \ \displaystyle\frac{\min_{k \in \mathcal{Q}}\lambda_{m+1}(k,x) -\max_{k\in \mathcal{Q}}\lambda_{m}(k,x)}{\tfrac{1}{2}(\min_{k\in \mathcal{Q}}\lambda_{m+1}(k,x)+\max_{k\in \mathcal{Q}}\lambda_{m}(k,x))}& \\ \\
     & \mbox{ s.t. } &  \bm{A}(k,x) u_j (k,x)= \lambda_j(k,x) \bm{M} u_j(k,x), \quad j = m, m+1, \ \mathrm{for~all~} k \in {\cal Q} \ .
\end{array}
\end{equation} The decision variables of problem $P$ are $x \in S \subset \mathbb{R}^{n_x}$, which represent the discretized material property of the design domain, as shown in the left of Figure \ref{fig_BGscheme}. The objective of problem $P$ is the eigenvalue {\em gap-midgap} ratio, which is defined as the difference between two prescribed consecutive eigenvalues divided by their mean (for scale invariance).  The constraints of problem $P$ are described in the equation line \eqref{eqBGopt}, which is shorthand for ``for each $k$ in the index set $\cal Q$, $\lambda_j(k,x)$ is the $j^{\mathrm{th}}$ ordered (generalized) eigenvalue (for $j=m,m+1$) of $\bm{A}(k,x)$ with respect to $\bm{M}$, for the given design variable $x$.''  Here $\cal Q$ is a particular governing index set, and is typically discretized to have the $n_k$ values ${\cal Q} =\{k_1, \ldots, k_{n_k}\}$, as shown in the right of Figure \ref{fig_BGscheme}.  (Indeed, in most bandgap problems $\cal Q$ indexes a discretization of the wave vectors which lie on the boundary of the Brillouin zone.)\medskip

Much of the details and derivation of the generalized eigensystem equation of \eqref{eqBGopt} are beyond the scope of the present paper.  However, we call out certain properties of $\bm{M}$ and the family of matrices $\bm{A}(k,x)$, as they are used in subsequent reformulations.
\begin{proposition}\label{proposition_EigenMat} Let $\lambda_1(k,x)\le \lambda_2(k,x) \le \ldots \le \lambda_{N}(k,x)$ denote the eigenvalues of the generalized eigensystem equation of \eqref{eqBGopt}, with corresponding normalized eigenfunctions $u_1(k,x), u_2(k,x), \ldots, u_N(k,x)$. Then $\bm{M}$ and $\bm{A}(k,x)$ have the following properties:
\begin{itemize}
\item[(i)] $\bm{M} \succ 0$,
\item[(ii)] $\bm{A}(k,x) \succeq 0,\ \mbox{for all } k\in \mathcal{Q},\ x\in S$,
\item[(iii)] $\bm{A}(k,x) = \bm{A}_0(k) + \sum_{i=1}^{n_x} \bm{A}_i(k)\ x_i,\ \mathrm{for~ all~} k\in \mathcal{Q}$,
\item[(iv)] $\lambda_i(k,x) \ge 0, \ \mathrm{for~all~ }i=1,\ldots,N$, and
\item[(v)] $u_i(k,x)^T \bm{M} u_j(k,x) = \delta_{ij},\ \mathrm{for~all~}i,j=1,\ldots,N$. \Halmos
\end{itemize}
\end{proposition}
Item $(i)$ of Proposition \ref{proposition_EigenMat} states that the (mass) matrix $\bm{M}$ is positive definite, while item $(ii)$ states that $\bm{A}(k,x)$ is positive semidefinite for any feasible design $x$, and for every $k \in {\cal Q}$. Item $(iii)$ states that the matrix $\bm{A}(k,x)$ depends affinely on the design variables $x$, for every $k \in {\cal Q}$. Items $(iv)$ and $(v)$ state that the generalized eigenvalues are nonnegative and the generalized eigenfunctions are $\bm{M}$-orthogonal, which are direct consequences of the previous three items.\medskip

Now let $\hat x \in S$ be given.  To ease the notation burden, we identify the finite set $\mathcal{Q} := \{k_1, \ldots, k_{n_k}\}$ with the counter $t \in \{1, \ldots, n_k\}$.  In \cite{men2010bandgap} it is shown how to construct operators ${\cal A}_{\ell,t}^{\hat{{x}}}(x) = A_{\ell,t,0}^{\hat{x}} + \sum_{i=1}^{n_x} A_{\ell,t,i}^{\hat{x}} x_i$ and ${\cal A}_{u,t}^{\hat{{x}}}( x) := A_{u,t,0}^{\hat{x}} + \sum_{i=1}^{n_x} A_{u,t,i}^{\hat{x}} x_i $, and also corresponding mass matrices $M_{\ell,t}^{\hat{{x}}}$ and $M_{u,t}^{\hat{{x}}}$, for each $t\in\{1, \ldots, n_k\}$, all of whose data depends on the current point $\hat x$, which are used to construct the following (convex) linear fractional semidefinite optimization problem (SDP):
\begin{equation}\label{eqSDP}
\begin{array}{lcll}
P_{SDP}^{\hat{x}}: &  \underset{x\in S,\lambda_{\ell},\lambda_u}{\max} & 2\displaystyle\frac{\lambda_u - \lambda_{\ell}}{\lambda_u+\lambda_{\ell}}& \\ \\
     & \mbox{ s.t. } &  {\cal A}_{\ell,t}^{\hat{{x}}}(x) \preceq \lambda_{\ell} M_{\ell,t}^{\hat{{x}}} \ , & t=1, \ldots, n_k \ , \\ [1.5ex]
     &               &  {\cal A}_{u,t}^{\hat{{x}}}(x) \succeq \lambda_u M_{u,t}^{\hat{{x}}} \ , & t=1, \ldots, n_k \ , \\ [1.5ex]
     &  &  \lambda_{\ell} \ge 0, \ \lambda_u \ge 0 \ .  &
\end{array}
\end{equation}
\begin{proposition}\label{proposition_FormSDP}
For a given $\hat{x}\in S$, the nonlinear nonconvex problem $P$ is locally approximated as the (convex) linear fractional semidefinite program $P_{SDP}^{\hat{x}}$ of \eqref{eqSDP}. \Halmos
\end{proposition}

Without going into the fine details, we note that $\lambda_{\ell}$ and $\lambda_u$ in \eqref{eqSDP} are intended to model $\max_{k\in \mathcal{Q}}\lambda_{m}(k,x)$ and $\min_{k\in \mathcal{Q}}\lambda_{m+1}(k,x)$ in \eqref{eqBGopt}, respectively, and that the two pairs of semidefinite inclusions in \eqref{eqSDP} locally model the $m^{\mathrm{th}}$ and $(m+1)^{\mathrm{st}}$ eigenvalue position for each $k\in \mathcal{Q} = \{k_1,\ldots,k_{n_k} \}$.  \medskip

Our goal herein is to apply the fabrication adaptivity paradigm to bandgap optimization problems.  Note that the objective function of \eqref{eqSDP} is at least as challenging as the largest eigenvalue function of Section \ref{eg1}.  Recall from the discussion in Section \ref{eg1} that the fabrication-adaptive counterpart function of the largest eigenvalue function is typically not computationally tractable.   We therefore proceed by replacing the semidefinite inclusions in \eqref{eqSDP} with linear inequality approximations, the methodology for which is described in Appendix \ref{subsec_PLrelax}, which yields the data $B^{\hat{{x}}}$, $C^{\hat{{x}}}$, $g^{\hat{{x}}}$, and $h^{\hat{{x}}}$ for the linear fractional optimization problem:
\begin{equation}\label{eq_lfp}
\begin{array}{lclr}
P_{LFP}^{\hat x}: &  \underset{ x\in S,\lambda_{\ell},\lambda_u}{\max} & 2\displaystyle \frac{\lambda_u - \lambda_{\ell}}{\lambda_u+\lambda_{\ell}}& \\ \\
     & \mbox{ s.t. } &  B^{\hat{{x}}}{x} + g^{\hat{{x}}}\le {e}\lambda_{\ell}, & \\ [1.5ex]
     &               &  C^{\hat{{x}}}{x} + h^{\hat{{x}}}\ge {e} \lambda_u, & \\ [1.5ex]
     &  &  \lambda_{\ell} \ge 0, \ \lambda_u \ge 0 \ .
\end{array}
\end{equation}
\begin{proposition}\label{proposition_FormLP}
For a given $\hat{x}\in S$, the linear fractional semidefinite program $P_{SDP}^{\hat{x}}$ of \eqref{eqSDP} is approximated as the linear fractional optimization problem $P_{LFP}^{\hat x}$ of \eqref{eq_lfp}. \Halmos
\end{proposition}

\noindent Here the two groups of semidefinite inclusions in \eqref{eqSDP} are replaced by $\mathcal{N}_B$ and $\mathcal{N}_C$ linear inequalities, respectively. 
The detailed methodology and derivation of the approximating linear inequalities are not the focus of the current work, but are nevertheless presented in Appendix \ref{subsec_PLrelax} for completeness.  In addition, the computational results presented in Appendix \ref{subsec_finality} show that solutions of \eqref{eq_lfp} are nearly as good and often are superior to those of \eqref{eqSDP}.  Our methodology is similar in spirit to that of \cite{sherali2002enhancing}, who replace semidefinite inclusions with linear inequalities to solve non-convex quadratic optimization problems on the simplex.\medskip

\subsection{Fabrication Adaptivity Formulation}
Notice that the linear fractional formulation in \eqref{eq_lfp} can be equivalently written in the following format which emphasizes the special piecewise linear structure of the objective function:
\begin{equation}\label{eq_lfpalt}
\begin{array}{clc}
&  \underset{x\in{S}}{\max} & f^{\hat{{x}}}({x}) :=  2\displaystyle \frac{\min\limits_{i\in \mathcal{I}}(C^{\hat{{x}}}{x} + h^{\hat{{x}}})_i - \max\limits_{j\in \mathcal{J}} (B^{\hat{{x}}}{x} + g^{\hat{{x}}})_j }{\min\limits_{i\in \mathcal{I}}(C^{\hat{{x}}}{x} + h^{\hat{{x}}})_i + \max\limits_{j\in \mathcal{J}} (B^{\hat{{x}}}{x} + g^{\hat{{x}}})_j},  \\ \\
\end{array}
\end{equation}
where ${\cal I} = \{1, \ldots, {\cal N}_C\}$ and ${\cal J} = \{1, \ldots, {\cal N}_B\}$.  Note in \eqref{eq_lfpalt} that the only constraint is the feasibility inclusion $x \in S$.  This special piecewise linear fractional function is of the exact structure as the problem discussed in Section \ref{eg5} (with the equivalence that now the objective is maximization rather than minimization and hence the roles of the ``max'' and ``min'' are switched in the fractional objective function).  Furthermore, it will also be the case in the bandgap application that the suppositions of Section \ref{eg5} are also satisfied, namely $S$ is compact and is conveyed as a system of linear inequalities $S=\{x : Ax \le b\}$, the norm $\|\cdot\|$ is a weighted $1$-norm $\|x\|:=\sum_{i=1}^{n_x} w_i|x_i|$ for positive weights $(w_1, \ldots, w_{n_x})$, and $(C^{\hat{{x}}}{x} + h^{\hat{{x}}})_i >0$ and $(B^{\hat{{x}}}{x} + g^{\hat{{x}}})_j > 0$ for all $x \in S$ and $i \in {\cal I}$ and $j \in \cal J$, respectively.  Therefore we can invoke the results in Section \ref{eg5} regarding computation and optimization of the fabrication-adaptive counterpart function $\tilde f^{\hat{{x}}}({x})$.  Let us see how this can be done.  From  Proposition \ref{separation} (with ``min'' and ``max'' appropriately interchanged) we have that:
\begin{equation}
\begin{array}{ll}
{f}^{\hat x}(x) & =  \min\limits_{{i \in \cal I},{j \in \cal J}}\  2\displaystyle\frac{(C^{\hat{x}}_i - B^{\hat{x}}_j) x + (h^{\hat{x}}_i - g^{\hat{x}}_j) }{ (C^{\hat{x}}_i + B^{\hat{x}}_j) x + (h^{\hat{x}}_i + g^{\hat{x}}_j)}\ .
\end{array}
\end{equation}

\noindent Furthermore, from \eqref{desk} and \eqref{fijast} we have that:
\begin{equation}
\begin{array}{ll}
\tilde{f}^{\hat x}(x) & =  \min\limits_{{i \in \cal I},{j \in \cal J}}\  \tilde{f}^{\hat{x}}_{i,j} (x)\ ,
\end{array}
\end{equation}
where for each $(i,j)\in \mathcal{I}\times \mathcal{J}$ we have:
\begin{equation}\label{ipod}
\begin{array}{llc}
\tilde{f}^{\hat{x}}_{ij}(x) := &\min\limits_{y } & \displaystyle  2\frac{(C^{\hat{x}}_i - B^{\hat{x}}_j) y + (h^{\hat{x}}_i - g^{\hat{x}}_j) }{ (C^{\hat{x}}_i + B^{\hat{x}}_j) y + (h^{\hat{x}}_i + g^{\hat{x}}_j)} \\[2ex]
& \mbox{s.t.} & \|y-x\| \le \delta \\[1ex]
& & y \in S\ .
\end{array}
\end{equation}
It also follows that $\tilde{f}^{\hat{x}}_{ij}(x)$ and $\nabla \tilde{f}^{\hat{x}}_{ij}(x)$ are computable via the convex optimization problems described in Propositions \ref{97} and \ref{supergradient3}, respectively, and that Algorithm FA of Table \ref{tb_PL} can be adapted (taking into account that the roles of min and max are switched and that the linear fractional pieces are indexed by pairs $(i,j)$) to solve the fabrication-adaptive optimization problem.   The format for  \eqref{eq_PL_FR0} becomes:
\begin{equation}\label{phone}\begin{array}{llrcl} \max\limits_{x\in S} \tilde f^{\hat x}(x) \ \ \ = \ \ \
& \underset{x,t}{\max} & t\\
& {s.t.} & \tilde{f}^{\hat x}_{i,j}(x) &\ge& t \ , \qquad \mathrm{for~all~} (i,j) \in {\cal I} \times {\cal J}\\
&& x &\in& S \ ,
\end{array}
\end{equation} and the linearization of \eqref{phone} at $\hat x$ then is:
\begin{equation}\label{fax}\begin{array}{llrcl}  \ \ \  \ \ \
& \underset{x,t}{\max} & t\\
& {s.t.} & \tilde{f}^{\hat x}_{i,j}(\hat x) + (\nabla\tilde{f}^{\hat x}_{i,j}(\hat x))^T(x-\hat x) &\ge& t \ , \qquad \mathrm{for~all~} (i,j) \in {\cal I} \times {\cal J}\\
&& x &\in& S \ .
\end{array}
\end{equation}  Table \ref{tb_PL2} presents the version of Algorithm FA, which we call Algorithm FA-B, for solving bandgap problems.

\begin{table}
\caption{fabrication-adaptive Optimization Algorithm for Bandgap Problems.\label{tb_PL2}}
\begin{center}
\begin{tabular}{|ll|}
\hline
& \textbf{Algorithm FA-B for Bandgap Problems}\\\hline
{\bf Step 1.} & Start with initial guess $\hat x := x^0$ and tolerance $\epsilon_{\mbox{tol}}$\\  & \\
{\bf Step 2a.}& Construct the linear operators for \eqref{eqSDP} based on $\hat x$, for  each $t\in\{1, \ldots, n_k\}$: \\
              & \quad ${\cal A}_{\ell,t}^{\hat{{x}}}(x) := A_{\ell,t,0}^{\hat{x}} + \sum_{i=1}^{n_x} A_{\ell,t,i}^{\hat{x}} x_i$\\
              & \quad ${\cal A}_{u,t}^{\hat{{x}}}( x) := A_{u,t,0}^{\hat{x}} + \sum_{i=1}^{n_x} A_{u,t,i}^{\hat{x}} x_i $\\
              & \quad $M_{\ell,t}^{\hat{{x}}}$ and $M_{u,t}^{\hat{{x}}}$\\
{\bf Step 2b.}& Construct the data for \eqref{eq_lfp} based on $\hat x$ and the linear operators from Step 2a:\\
              & \quad $B^{\hat{{x}}}$, $g^{\hat{x}}$, $C^{\hat{{x}}}$, and $h^{\hat{x}}$\\
{\bf Step 3.} & For each $(i,j) \in {\cal I} \times {\cal J}$, do:\\
              & \quad Compute the function value $\tilde{f}_{i,j}(\hat{x})$ of \eqref{ipod} via \eqref{eq_PLi}\\
              & \quad Compute first-order information $\nabla\tilde{f}_{i,j}(\cdot)$ at $\hat x$ via \eqref{eq_PLi_grad}\\ & \\
{\bf Step 4.} & Form the linear optimization problem \eqref{fax} \\ & \\
{\bf Step 5.} & Solve \eqref{fax} for an optimal solution $x^*$\\ & \\
{\bf Step 6.} & If $\|x^* - \hat x\| \le \epsilon_{\mbox{tol}}$, stop.\\
              & Else update $\hat x \leftarrow x^*$ and go to {\bf Step 2.}\\
\hline
\end{tabular}
\end{center}
\end{table}

\section{Computational Results for Fabrication-Adaptive Optimization}\label{pcd}
\label{sec_res}

We first test the intended effectiveness of Algorithm FA (Table \ref{tb_PL}) on random problems.  These results are presented in Section \ref{eaves}.  We then apply Algorithm FA-B to a variety of bandgap problems that arise in photonic crystal design, which was the problem class that engendered this line of research.  We show via several examples how the solutions produced by Algorithm FA-B succeed in producing improved fabricable design solutions compared to solutions based on the optimal solution of the original design optimization problem.  These results are presented in Section \ref{saigal}.

\subsection{Computation on Random Piecewise Linear Fractional Problems}\label{eaves}
We tested Algorithm FA on randomly generated instances of the specially structured piecewise linear fractional optimization problem \eqref{sandy} presented in Section \ref{eg5}.  All problems were generated using $n=50$, $n_{\cal I} = 20$, and $n_{\cal J} = 30$, and $S$ set to a unit hypercube of the form $S:=\{x \in \mathbb{R}^n : 1.0 \le x_i \le 2.0, \ i=1, \ldots, n \}$.  All of the components of $ a^i, c^j, g_i, h_j$ were chosen randomly from the uniform distribution ${\cal U}[0,1]$, for $i = 1, \ldots, n_{\cal I}$, $j = 1,\ldots, n_{\cal J}$.  The choice of $S$, the range of the data and the size of the generated problems are in close agreement with problems arising in photonic crystal design problems to be discussed in Section \ref{saigal}.\medskip

Given a randomly generated instance of \eqref{sandy}, let us denote the optimal solution of the original optimization problem \eqref{sandy} as:
\begin{equation}
\begin{array}{l}
x^{\ast}_O := \arg\ \min\limits_{x\in S} f(x) \ .
\end{array}
\label{eq_xo}
\end{equation}
To construct the fabrication-adaptive optimization problem, we used the norm $\|\cdot\|:= \|\cdot\|_1$ and set $\delta = 5.0$, and note that this choice allows for changing up to $10\%$ ($5/50$) of the components $x_i$ of $x$ by one unit, which is the full range of $x_i$ by definition of $S$.  Let $x^{\ast}_{FA}$ denote the computed solution of the fabrication-adaptive optimization problem using Algorithm FA.\medskip

Note that $x^{\ast}_O$ is the optimal solution of \eqref{sandy} which is solvable as a linear optimization problem.  In contrast, $x^{\ast}_{FA}$ is not necessarily the (global) optimal solution of FA optimization problem \eqref{eq_PL_FR0} and is computed using Algorithm FA (Table \ref{tb_PL}).  Here and in what follows we used the Gurobi Optimizer \citep{gurobi} to solve all linear optimization problems.\medskip

\begin{figure}[htdp]
\begin{center}
\includegraphics[scale=0.5]{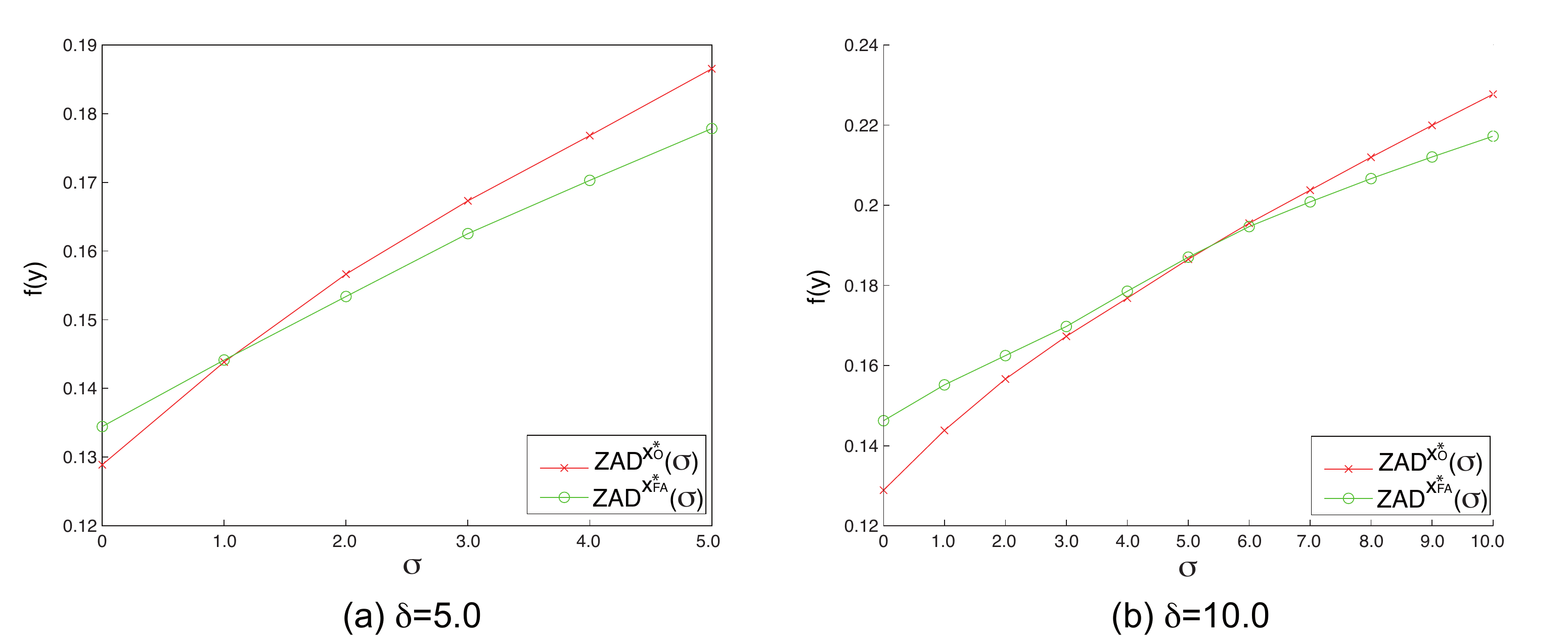}
\caption{Adaptivity of $x^{\ast}_{O}$ and $x^{\ast}_{FA}$ to nearby \emph{adversarial} solutions, as defined in \eqref{xadv}.}
\label{fig_eg1}
\end{center}
\end{figure}

We generated $20$ random instances of \eqref{sandy}. Since the results of all $20$ instances lead to the same conclusion, we only discuss detailed results from a particular one.  For this problem instance the optimal objective function value of the original problem is $f(x^{\ast}_O) = 0.129$, whereas the value of the original objective function evaluated at the fabrication-adaptive solution is $f(x^{\ast}_{FA})=0.135$, which is inferior (larger) to that of the optimal value as expected. We tested the adaptivity of $x^{\ast}_{O}$ and $x^{\ast}_{FA}$ as follows.  For these two solutions under consideration, and a given value of $\sigma \in [0,\delta]$, we compute the most conservative objective function value $f(y)$ among all solutions $y$ for which $\|y-\hat x\| \le \sigma$ and $y \in S$.  That is, we compute:
\begin{equation}\label{xadv}
\begin{array}{ll}
\mathrm{ZAD}^{\hat x}(\sigma) := & \max\limits_{{y}} \ \ f({y})   \\
  &  \mbox{s.t.}\ \ \ \ \|{y}-\hat x\|\leq \sigma  \\
  & \qquad \ \ {y}\in S \ ,
\end{array}
\end{equation}\medskip

\noindent for $\hat x = x^{\ast}_{O}$ and $\hat x = x^{\ast}_{FA}$.  One can interpret \eqref{xadv} as computing the worst solution $y$ whose distance from $\hat x$ is at most $\sigma$.  In the absence of an intelligent method for adapting a solution $\hat x$, \eqref{xadv} essentially assumes the solution $\hat x$ will be adapted to a nearby solution $y$ in an {\em adversarial} manner (hence the choice of notation ``$ZAD$'' in \eqref{xadv}).  The values of $\mathrm{ZAD}^{\hat x}(\sigma)$ were computed using Steps 2 -- 4 of Algorithm FA in Table \ref{tb_PL}.  Plots of $\mathrm{ZAD}^{x^{\ast}_{O}}(\sigma)$ and $\mathrm{ZAD}^{x^{\ast}_{AF}}(\sigma)$ for $\sigma \in [0,\delta]$ are shown in Figure \ref{fig_eg1}$(a)$.  For small values of $\sigma$, the range of adversarial solutions is small, and hence the superior original objective function value of $x^{\ast}_{O}$ yields $\mathrm{ZAD}^{x^{\ast}_{O}}(\sigma) < \mathrm{ZAD}^{x^{\ast}_{FA}}(\sigma)$.  However, as the values of $\sigma$ increases, the superior adaptability of the solution $x^{\ast}_{FA}$ is revealed.  When the range of adversarial solutions is larger, $\mathrm{ZAD}^{x^{\ast}_{O}}(\sigma) > \mathrm{ZAD}^{x^{\ast}_{FA}}(\sigma)$, showing that nearby adversarial solutions of $x^{\ast}_{FA}$ are superior to those of $x^{\ast}_{O}$.  These plots reveal that the solution $x^{\ast}_{FA}$ is indeed effectively more adaptive, with the advantage growing as the allowable range of nearby solutions grows.  We also repeated this computational exercise using $\delta = 10.0$.  The resulting plots of $\mathrm{ZAD}^{x^{\ast}_{O}}(\sigma)$ and $\mathrm{ZAD}^{x^{\ast}_{AF}}(\sigma)$ for the case of $\delta = 10.0$ are shown in Figure \ref{fig_eg1}$(b)$.  Notice that the results for the case $\delta = 10.0$ further reinforce the above observations.\medskip

Of the $20$ randomly generated instances of \eqref{sandy}, all exhibited similar effectiveness of the adaptivity of $x^{\ast}_{FA}$ in terms of $\mathrm{ZAD}^{x^{\ast}_{O}}(\sigma)$ and $\mathrm{ZAD}^{x^{\ast}_{FA}}(\sigma)$.

\subsection{Computational Experience on Bandgap Problems in Photonic Crystal Design}\label{saigal}

In this section we present results from applying the fabrication-adaptive optimization model to bandgap problems in photonic crystal design, which was the originator of our need to pursue this line of research.  As briefly reviewed in Section \ref{sec_FA-BOP}, the goal is to optimize the bandgap between two consecutive eigenvalues, where the bandgap is the largest gap that separates the two eigenvalues over all values of $k$ in the governing index set $\cal Q$.  We seek to solve fabrication-adaptive models \eqref{phone} using Algorithm FA-B (Table \ref{tb_PL2}).  There are many different types of bandgap optimization problems that one can construct as well as different schemes themselves for constructing bandgap optimization problems in photonic crystal design.  For example, one arrives at different bandgap optimization problems depending on which eigenvalue gap one seeks to optimize (the $m^{\mathrm{th}}$ bandgap, defined as the relative gap between the $m^{\mathrm{th}}$ and the $(m+1)^{\mathrm{st}}$ eigenvalues for $m=1, \ldots, N-1$), the choice of polarization (TE or TM or complete (TEM) polarization), and the  lattice structure of the photonic crystal (typically a square lattice or a triangular lattice).  Different combinations of these choices lead to different bandgap optimization problems with different optimal solutions. Among the numerous topological varieties in the optimal structures derived from different bandgap problems, one often encounters solutions that are either not fabricable or pose onerous fabrication challenges due to thin connectors, small features, rough edges, isolated structures, and other related solution configurations.\medskip

Among the roughly $60+$ bandgap problems that we have solved, the original optimal solutions of at least 15 are not fabricable without post-processing modification.  Furthermore, most of the non-fabricable solutions are for optimization problems for complete (TEM) bandgaps or other multiple-bandgap problems \citep{men2011thesis, men2011design}.  We applied the fabrication adaptive optimization paradigm and algorithms to most of these problems.  Herein we report on some of our computational experience to address questions such as: (i) how sensitive are solutions of the original problem to fabrication adaptivity modifications?, (ii) how good are the solutions computed when solving the fabrication-adaptive optimization problem?, and (iii) how do fabrication adapative solutions compare to solutions of the original optimization problem?\medskip

For a given problem instance, let $x_O^{\ast}$ be the optimal solution of the original bandgap problem \eqref{eqSDP} (or \eqref{eq_lfp}) and let $x^{\ast}_{FA}$ be the computed solution of the fabrication-adaptive optimization problem using Algorithm FA-B.  In the case when $x_O^{\ast}$ and/or $x^{\ast}_{FA}$ are not fabricable, we applied manual changes of these solutions to produce fabricable solutions $y_{FA}$ and/or $y_O$ by using our own problem-domain common sense to modify pixel values to create more fabricable designs.  The manual changes we employed are of the ordinary variety such as removing thin rods, removing small features, smoothing boundaries of materials, and straightening inner edges of material boundaries.  All of these modifications can be easily accomplished with standard image processing filters.\medskip

\begin{figure}[htdp]
\begin{center}
\includegraphics[scale=0.5]{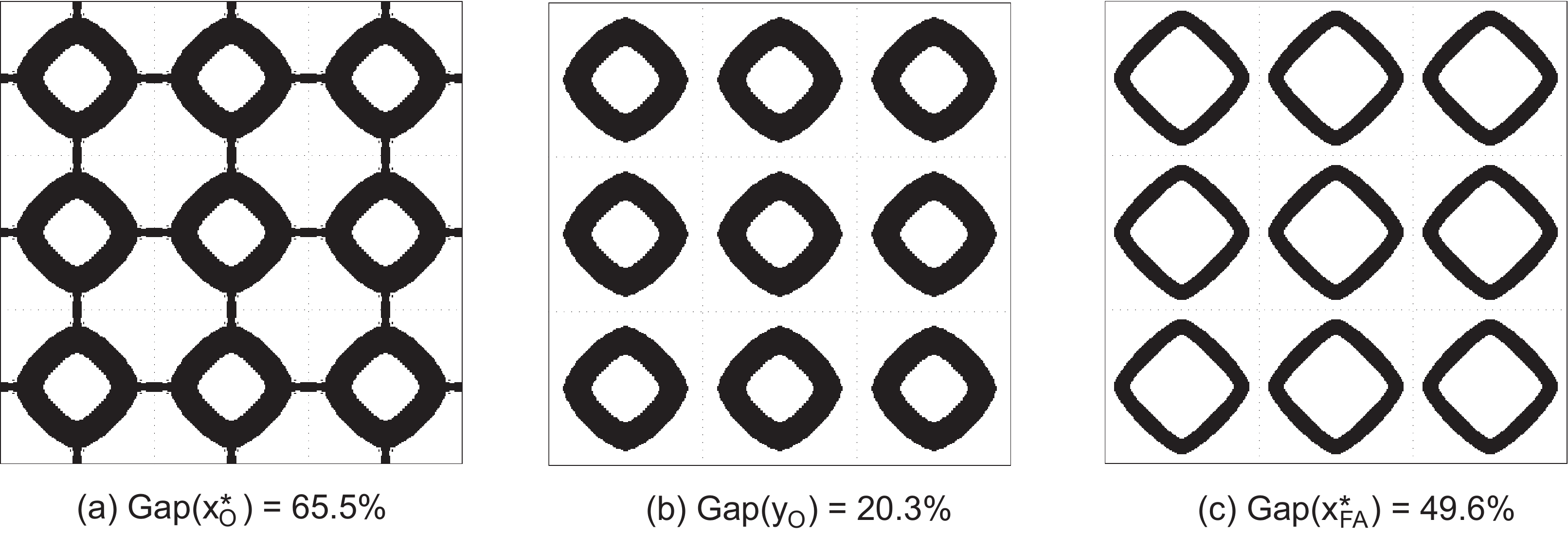}
\end{center}
\caption{Solutions to the $2$nd TE bandgap problem in the square lattice. (a) is the original optimal design $x^{\ast}_O$, with Gap$(x^{\ast}_O) =  65.5\% $;  (b) is the solution $y_O$ which is a manual modification of $x^{\ast}_O$, with Gap$(y_O) =  20.3\% $, and $3\%$ of pixels being modified; (c) is the computed solution $x^{\ast}_{FA}$ of using Algorithm FA-B, with Gap$(x^{\ast}_{FA}) =  49.6\%$, using $\delta = 3\%$.}
\label{fig_sq_te2}
\end{figure}

As in Section \ref{eaves}, let $x_O^{\ast}$ be the optimal solution of the original bandgap problem \eqref{eqSDP} (or \eqref{eq_lfp}).  An example of the poor performance of $x_O^{\ast}$ after modification for fabrication arises in solving for the $2$nd TE bandgap in the square lattice, as shown in Figure \ref{fig_sq_te2}.  Figure \ref{fig_sq_te2}$(a)$ shows the design $x_O^{\ast}$ which has a bandgap of $65.6\%$, but which contains very thin rods that are challenging to fabricate.  One can remove the thin rods by modifying $3\%$ of the pixels, yielding the design $y_O$ shown in part $(b)$ of Figure \ref{fig_sq_te2}.  However, this small modification of the fabricable solution $y_O$ drastically reduces the bandgap, from $65.6\%$ down to $20.3\%$.  Figure \ref{fig_sq_te2}$(c)$ shows the computed solution $x^{\ast}_{FA}$ of the fabrication-adaptive optimization problem using Algorithm FA-B, for the value of  $\delta=3\%$. Not only is the bandgap for $x^{\ast}_{FA}$ much higher ($49.6\%$) than that of $y_O$, but it actually is fabricable as is. This is very fortunate, but perhaps ``accidental'', as it is more typical that the solution $x^{\ast}_{FA}$ would need to be modified to a nearby solution as we will see below.  The main point of this example is to show how unadaptable the solution $x_O^{\ast}$ can be to modification that will make it fabricable without unduly reducing the size of the bandgap.\medskip

\begin{figure}[htdp]
\begin{center}
\includegraphics[scale=0.5]{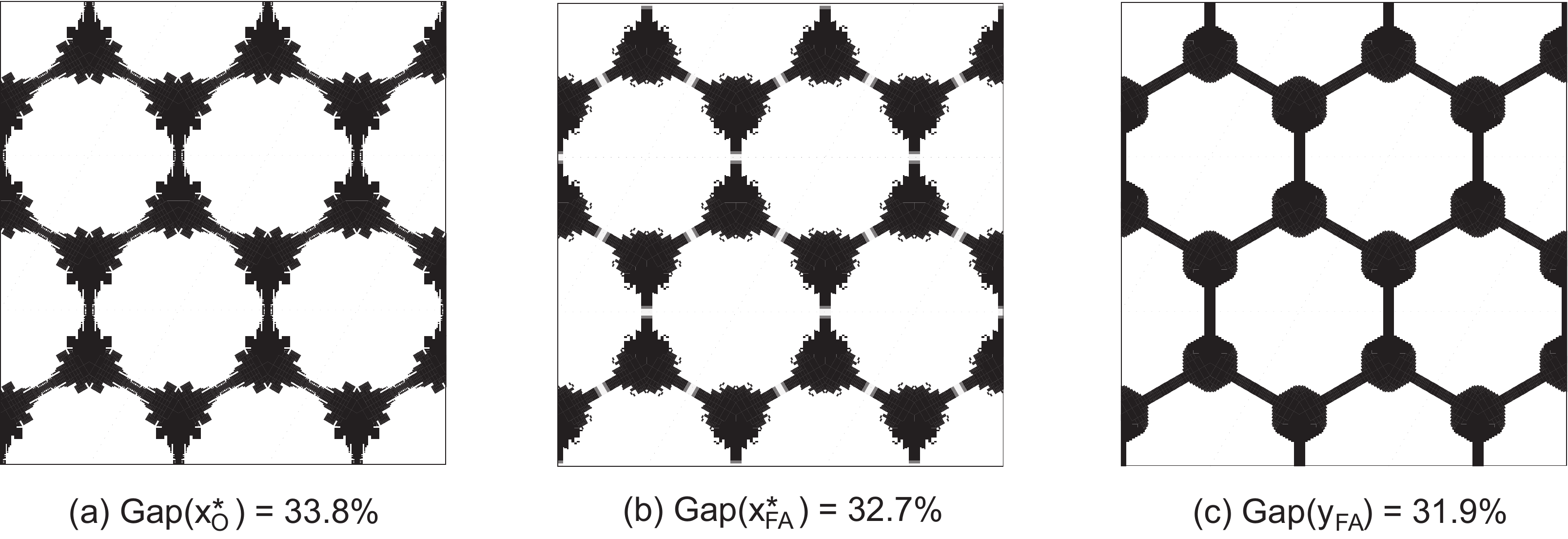}
\end{center}
\caption{Solutions to the complete (TEM) bandgap problem involving the $1$st TE bandgap and the $2$nd TM bandgap, in the triangular lattice. (a) is the original optimal design $x^{\ast}_O$, with Gap$(x^{\ast}_O) =  33.8\% $;  (b) is the computed solution $x^{\ast}_{FA}$ using Algorithm FA-B, with Gap$(x^{\ast}_{FA}) =  32.7\%$, using $\delta = 5\%$; (c) is the solution $y_{FA}$ which is a manual modification of $x^{\ast}_{FA}$, with Gap$(y_{FA}) =  31.9\%$, and $4\%$ modification.}
\label{fig_tr_te1tm2}
\end{figure}

An example of the quality of solutions computed using Algorithm FA-B is shown in Figure \ref{fig_tr_te1tm2}, which shows solutions to the complete (TEM) bandgap problem involving the $1$st TE bandgap and the $2$nd TM bandgap, in the triangular lattice.  For this bandgap problem the bandgap for the original (non-fabricable) solution $x^*_O$ (Figure \ref{fig_tr_te1tm2}$(a)$) is $33.3\%$.  Figure \ref{fig_tr_te1tm2}$(b)$ shows the computed solution $x^{\ast}_{FA}$ of the fabrication-adaptive optimization problem using Algorithm FA-B, for the value of  $\delta=5\%$.  Figure \ref{fig_tr_te1tm2}$(c)$ shows the manually modified solution $y_{FA}$ of $x^{\ast}_{FA}$, with about $4\%$ of pixels being modified. $y_{FA}$ is significantly more fabricable than $x^*_O$, yet its bandgap value $31.9\%$ is only modestly decreased from that of the original (non-fabricable) solution.\medskip

\begin{figure}[htdp]
\begin{center}
\includegraphics[scale=0.5]{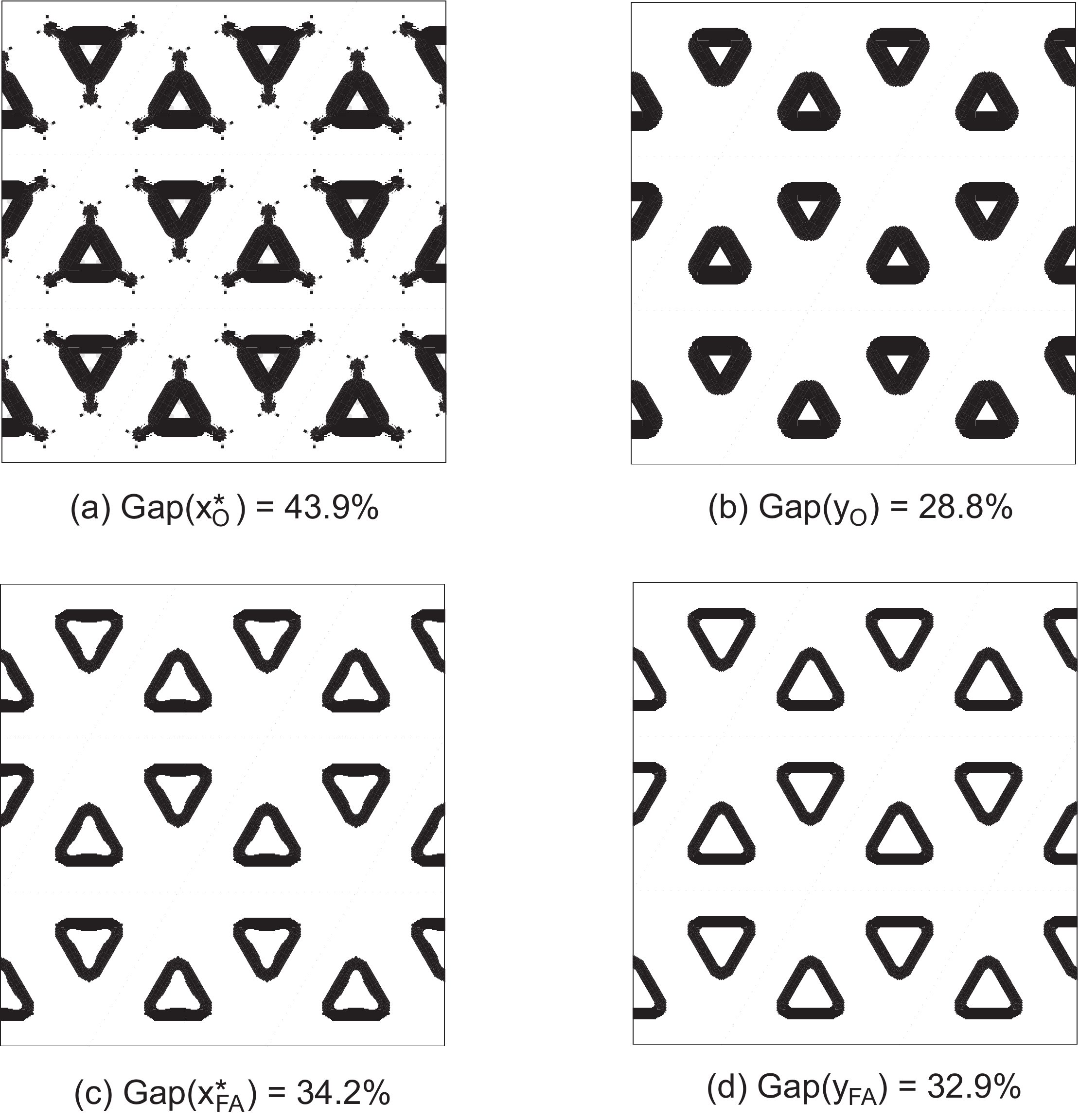}
\end{center}
\caption{Designs with $5^{\mathrm{th}}$ TE eigenbandgap in a triagnular lattice. (a) original optimal design $x^{\ast}_O$, with $ Gap(x^{\ast}_O) =  43.9\% $;  (b) manually modified design $y_O$ based on the original design, with $ Gap(y_O) =  28.8\% $, and $5\%$ modifications; (c) Fabrication-Adaptive optimal design $x^{\ast}_{FA}$, with $ Gap(x^{\ast}_{FA}) =  34.2\%$, and $\delta_{FA} = 5\%$; (d) manually modified design $y_{FA}$ based on the FA optimal design, with $ Gap(y_{FA}) =  32.9\% $, and $0.8\%$ modification.}
\label{fig_tr_te5}
\end{figure}

Our last two examples illustrate the comparative value of the fabrication-adaptive optimization approach. First we solve for the $5^{\mathrm{th}}$ TE bandgap in the triangular lattice, as shown in Figure \ref{fig_tr_te5}.  By simply eliminating the small features of the original optimal solution $x_O^{\ast}$ shown in Figure \ref{fig_tr_te5}$(a)$ (which comprise $5\%$ of the pixels), the bandgap of the manually modified solution $y_O$ is sharply decreases from $43.9\%$ to $28.8\%$.  However, the fabrication-adaptive computed solution using Algorithm FA-B (using the same modification allowance $\delta=5\%$ of pixels) yields the solution $x^{\ast}_{FA}$ (shown in Figure \ref{fig_tr_te5}$(c)$).  The manual modification of this solution is $y_{FA}$, and is shown in Figure \ref{fig_tr_te5}$(d)$.  The modified solution $y_{FA}$ is designed so that the inner edges of the triangular structures in $x^{\ast}_{FA}$ are straight, in order to make the resulting solution more fabricable.  The resulting bandgap of the modified solution $y_{FA}$ is $32.9\%$, which is better than that of the fabricable solution $y_O$ based on the original optimal solution.

\begin{figure}[htdp]
\begin{center}
\includegraphics[scale=0.5]{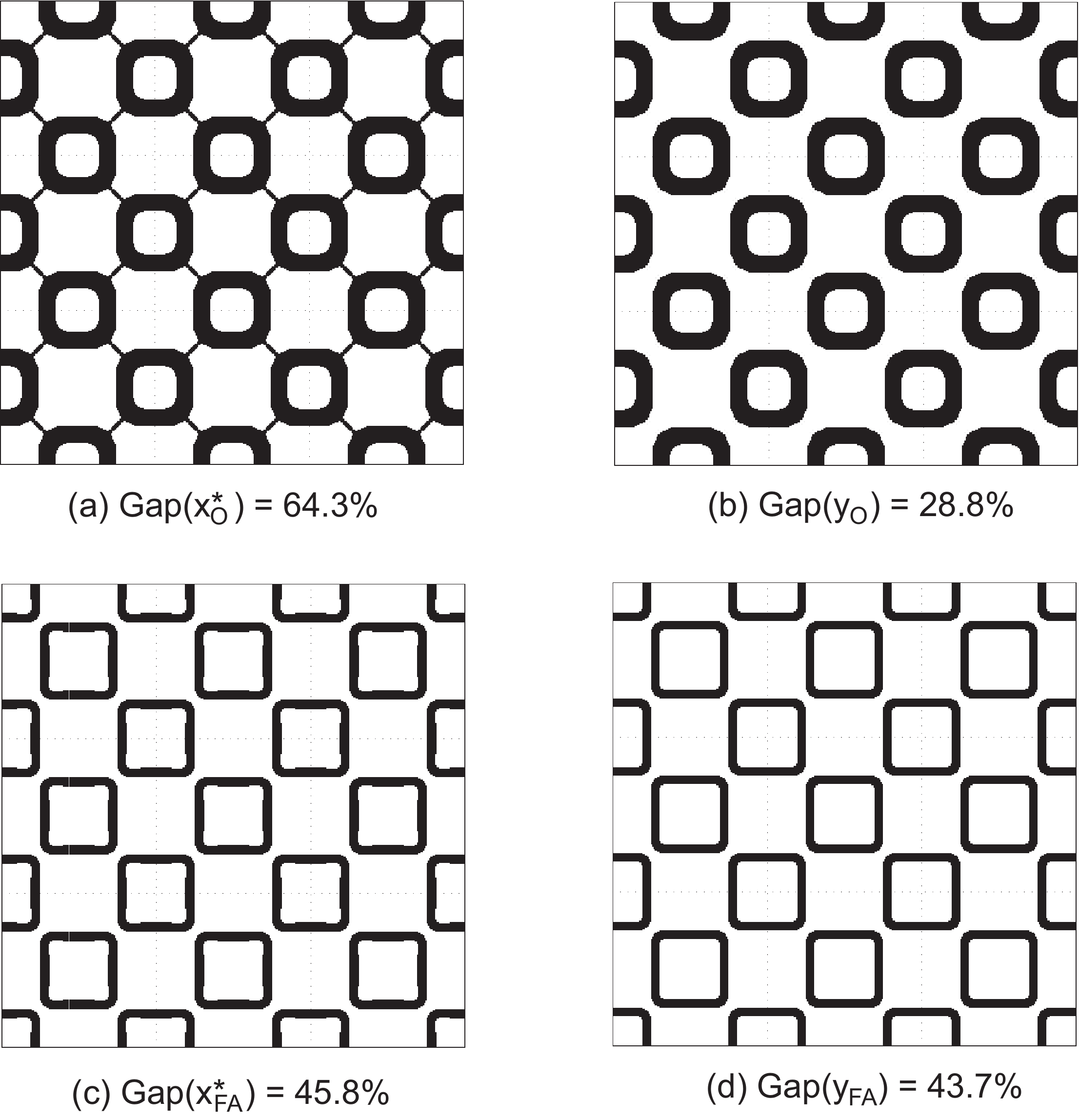}
\end{center}
\caption{Solutions to the $4^{\mathrm{th}}$ TE bandgap in the square lattice. (a) is the original optimal design $x^{\ast}_O$, with Gap$(x^{\ast}_O) =  64.3\% $; (b) is the solution $y_O$ which is a manual modification of $x^{\ast}_O$, with Gap$(y_O) =  28.8\% $, and $1.2\%$ of pixels being modified;  (c) is the computed solution $x^{\ast}_{FA}$ using Algorithm FA-B, with Gap$(x^{\ast}_{FA}) =  45.8\%$, using $\delta = 5\%$; (d) is the solution $y_{FA}$ which is a manual modification of $x^{\ast}_{FA}$, with Gap$(y_{FA}) =  43.7\% $,  and also $1.2\%$ of pixels being modified.}
\label{fig_sq_te4}
\end{figure}

A similar situation occurs to solutions of the $4^{\mathrm{th}}$ TE bandgap in the square lattice.  Figure \ref{fig_sq_te4}$(a)$ and Figure \ref{fig_sq_te4}$(b)$ show the original optimal solution and its manual modification to a fabricable solution by removing the thin rods, which reduces the original bandgap from $64.3\%$ to $28.8\%$.  In contrast, the computed solution $x^{\ast}_{FA}$ shown in Figure \ref{fig_sq_te4}$(c)$ is truly more adaptive to modification for fabrication.  Figure \ref{fig_sq_te4}$(d)$ shows the manually modified solution $y_{FA}$ of $x^{\ast}_{FA}$, which was done by straightening the inner sides of the square structures. In the manual modifications of both solutions $x^{\ast}_{O}$ and $x^{\ast}_{FA}$, the fraction of pixels modified was very small, both are roughly $1.2\%$.  However, in the case of the original solution, the bandgap was significantly reduced (from $64.3\%$ to $28.8\%$), while in the case of the fabrication-adaptive solution the bandgap reduction was minor  (from $45.8\%$ to $43.7\%$).\medskip

These examples show that while manual modification of solutions might appear to be minor, the negative effect on the bandgap can be very significant at solutions to the original problem, but (as intended) are less significant at computed solutions to the fabrication-adaptive optimization problem.\medskip

\section{Conclusions}

We have introduced the \emph{fabrication-adaptive} optimization modeling paradigm \eqref{eq_fr-cp}-\eqref{eq_fr_opt}, which stems from the {\em robust regularization} operation on functions \citep{lewis2002}.  The FA modeling paradigm does not necessarily yield a convex optimization model even when the original optimization problem is convex.  Hence, we examined a variety of special structures on functions, feasible regions, and norms, for which computation is tractable, and we developed an algorithmic scheme for solving certain FA optimization problems that arise from piecewise linear fractional optimization.  We first tested the FA paradigm and algorithm on randomly generated problems to show some general behavior of solutions.  We next applied our methodology to bandgap optimization problems in photonic crystal design, which were the originating class of problems that engendered this line of research. These bandgap problems were originally modeled using SDP formulations of iteration-specific approximation problems.  To apply the FA framework, we developed piecewise linear approximations of the semidefinite inclusions, which worked surprisingly well and enabled replacing SDP inclusions with linear inequalities that yielded linear optimization problems.  We used the FA model and algorithm to compute significantly improved fabricable designs of a variety of bandgap optimization problems in photonic crystal design.\medskip

As mentioned above, the success of piecewise linear approximations of the semidefinite inclusions in bandgap optimization models is counter to traditional notions that such approximations are crude at best.  It is unclear at this point whether the success of our simple LP/SDP approximation is due to the very specific structure of bandgap design problems and the resulting eigenvalue bound inclusions.  Future research on our agenda includes other applications and extensions of fabrication-adaptive optimization, as well as exploration aimed at understanding the possible reach of success of the LP/SDP approximation method described in the Appendix \ref{sec_PLrelax} of this paper.

%


%

%
%
%
\begin{APPENDICES}
\SingleSpacedXI

\section{Relation of Fabrication-Adaptive Optimization Model to Robust Optimization}\label{robustopt}

In the special case when $S = \mathbb{R}^n$, one can re-formulate the fabrication-adaptive optimization problem \eqref{eq_fr_opt} as a particular instance of a robust optimization problem \citep{ben2009robust}, at least conceptually.  Let us see how this can be done.  We first use the change of variable $d:= y-x$ and note that when $S = \mathbb{R}^n$ we can re-write \eqref{eq_fr-cp} as:
\begin{equation}\label{eq_fr-cp1}
\begin{array}{ll}
\tilde{f}({x}) = & \max_{{d}}\ \  f({x+d})   \\
  &  \mbox{s.t.}\ \ \ \ \|d\|\leq \delta \  \\
  & \qquad \ \ {d}\in \mathbb{R}^n \ .
\end{array}
\end{equation}
Considering $d$ as the ``data'' we can define the function $\hat f_d(x):=f(x+d)$ for $x,d \in \mathbb{R}^n$, where the ``data'' $d$ parametrically defines the function $\hat f_d(\cdot)$.  Then notice that the level set condition `` $\tilde f(x) \le t$ '' obeys: \begin{equation}\label{eq_fr-cp2}
\begin{array}{ll} \tilde f(x) \le t \ \ \ \ \ \ \Leftrightarrow \ \ \ \ \
\hat f_d(x) \le t & \mbox{for all~} d \in B(0,\delta) \ ,
\end{array}
\end{equation}
where $B(c,r)$ denotes the ball centered at $c$ with radius $r$.  Therefore, we can write the fabrication-adaptive optimization problem \eqref{eq_fr_opt} as:
\begin{equation}\label{eq_fr_opt2}
\begin{array}{lll}
\tilde{z}^{\ast} = & \min_{x,t} \ \  t  \\
  &  \mbox{s.t.} \ \  \hat f_d(x) \le t & \mbox{for all~} d \in B(0,\delta) \ .
\end{array}
\end{equation}
Observe that \eqref{eq_fr_opt2} corresponds exactly to a robust optimization model with uncertain ``data'' $d$ used as the ``data'' parameter of the function $\hat f_d(x):=f(x+d)$, and with the uncertainty set ${\cal U}:= B(0,\delta)$.  In the language of robust optimization, the constraints of \eqref{eq_fr_opt2} immunize the inequality ``$f(x) \le t$'' over all possible values of the data $d$ in the uncertainty set ${\cal U}:= B(0,\delta)$.  If $f(\cdot)$ is a convex function, then $\hat f_d(\cdot)$ is convex for any $d$, whereby $\tilde f(x) = \max_{d \in B(0,\delta)} f_d(x)$ is also convex as it is the pointwise maximum of convex functions.  (This also provides an alternate proof that the fabrication-adaptive optimization problem \eqref{eq_fr_opt} is a convex optimization problem when $S=\mathbb{R}^n$.)\medskip

When $S \neq \mathbb{R}^n$, we show that the above analysis breaks down. In the general case of $S \subset \mathbb{R}^n$ the fabrication-adaptive optimization problem \eqref{eq_fr_opt} can be re-written as:
\begin{equation}\label{eq_fr_opt3}
\begin{array}{lll}
\tilde{z}^{\ast} = & \min_{x,t}  \ \ t  \\[1ex]
  &  \mbox{s.t.} \ \ f_d(x) \le t & \mbox{for all~} d \in B(0,\delta)\cap\left(S -\{x\} \right) \\[1ex]
  & \qquad \ \ x \in S \ .
\end{array}
\end{equation}  Now notice in \eqref{eq_fr_opt3} that the corresponding ``uncertainty set'' is now ${\cal U} = {\cal U}(x):= B(0,\delta)\cap \left(S -\{x\}\right)$ which depends on the decision variable $x$.  The lack of independence of the uncertainty set ${\cal U} = {\cal U}(x)$ from the value of the variable $x$ leads to the potential for the problem \eqref{eq_fr_opt3} to be non-convex even when $f(\cdot)$ is convex and the feasible region $S$ is convex.  It was already shown in Section \ref{ci} that one can easily construct such an instance where the resulting fabrication-adaptive optimization problem is not convex and is not even quasi-convex. \medskip

\section{Relaxation and Reformulation of the bandgap Optimization Problem}
\label{sec_PLrelax}

\subsection{Approximation of the Eigenvalue Bounds using Linear Inequalities}\label{subsec_PLrelax}
While we would like to apply the fabrication-adaptive methodology to the bandgap optimization problem, the third example in Section \ref{eg1} illustrates the challenges in doing so.  The objective function of the bandgap optimization problem \eqref{eqSDP} is at least as complicated as the largest eigenvalue function \eqref{eq_eig}, whose fabrication-adaptive counterpart is not generally tractable to compute as discussed in Section \ref{eg1}.  However, Example \ref{eg5} shows that if $f(\cdot)$ is a special piecewise linear fractional function, then its fabrication-adaptive counterpart $\tilde f(\cdot)$ is tractable to compute.  We therefore propose to replace the eigenvalue bounds in \eqref{eqSDP}, which are modeled with semidefinite inclusions, with piecewise linear approximations that are modeled with linear inequalities, thereby replacing \eqref{eqSDP} with a linear fractional optimization problem of the form \eqref{eq_lfp}.   We carry out this step as follows.\medskip

The matrix data in \eqref{eqSDP} are constructed as ``reduced'' stiffness and mass matrices based on the current iterate $\hat x \in S$, and are given by  (see \cite{men2010bandgap}):
\begin{equation}\label{eqRedMat}
\begin{array}{ll}
\mathcal{A}_{\ell,t}^{\hat{{x}}}(x) &:= A_{\ell,t,0}^{\hat{x}} + \sum_{i=1}^{n_x} A_{\ell,t,i}^{\hat{x}} x_i := \Phi_{\ell}^{\hat{ x}} ({k}_t)^{\ast} \bm{A}_0({k}_t) \Phi_{\ell}^{\hat{ x}} ({k}_t) + \sum_{i=1}^{n_x}\bigl( \Phi_{\ell}^{\hat{ x}} ({k}_t)^{\ast} \bm{A}_i({k}_t) \Phi_{\ell}^{\hat{ x}} ({k}_t) \bigr) x_i,  \\[1ex]
M_{\ell,t}^{\hat{x}} &:=  \Phi_{\ell}^{\hat{ x}} ({k}_t)^{\ast} \bm{M} \Phi_{\ell}^{\hat{ x}} ({k}_t), \\[1ex]
\mathcal{A}_{u,t}^{\hat{{x}}}( x) &:= A_{u,t,0}^{\hat{x}} + \sum_{i=1}^{n_x} A_{u,t,i}^{\hat{x}} x_i := \Phi_{u}^{\hat{ x}} ({k}_t)^{\ast} \bm{A}_0({k}_t) \Phi_{u}^{\hat{ x}} ({k}_t) + \sum_{i=1}^{n_x}\bigl( \Phi_u^{\hat{ x}} ({k}_t)^{\ast} \bm{A}_i({k}_t) \Phi_u^{\hat{ x}} ({k}_t) \bigr) x_i ,\\[1ex]
M_{u,t}^{\hat{x}} &:=\Phi_u^{\hat{ x}} ({k}_t)^{\ast} \bm{M} \Phi_u^{\hat{ x}} ({k}_t) ,
\end{array}
\end{equation}
 for $t=1, \ldots, n_k$. The subspace matrices $\Phi_{\ell}^{\hat{ x}} ({k}_t) $ and $\Phi_u^{\hat{ x}} ({k}_t) $ consist columnwise of the ``important'' eigenfunctions,
  \begin{equation}
\Phi_{\ell}^{\hat{ x}} ({k}_t) = \left[ u_a(k_t,\hat{x}),\ldots,u_{m}(k_t,\hat{x})\right], \qquad
\Phi_u^{\hat{ x}} ({k}_t) = \left[ u_{m+1}(k_t,\hat{x}),\ldots,u_{b}(k_t,\hat{x})\right].
 \end{equation}
The following result is obtained as a consequence of the derivation of these matrices from the Finite Element Method and from the fact that the basis sets $\Phi_\ell^{\hat{ x}} ({k}_t) $, $\Phi_u^{\hat{ x}} ({k}_t) $ are $\bm{M}-$orthogonal bases, see also Proposition \ref{proposition_EigenMat}.
\begin{proposition}\label{theorem_RedMat}
For all $x\in S$, the reduced stiffness and mass matrices \eqref{eqRedMat} satisfy:
\begin{enumerate}
\item $\mathcal{A}_{\ell,t}^{\hat{x}}(x)\succeq 0,\mbox{and } \mathcal{A}_{u,t}^{\hat{x}}(x)\succeq 0$ for $t=1, \ldots, n_k$,
\item $\mathcal{A}_{\ell,t}^{\hat{x}}(x)\succ 0, \mbox{and } \mathcal{A}_{u,t}^{\hat{x}}(x)\succ 0$ for $k_t \neq 0$, and
\item $M_{\ell,t}^{\hat{x}}\succ 0, \mbox{and } M_{u,t}^{\hat{x}}\succ 0$ for $t=1, \ldots, n_k$.\Halmos
\end{enumerate}
\end{proposition}
We next note that the semidefinite inclusions in \eqref{eqSDP} can be rewritten as:
\begin{equation}\label{eq_sdp2lp1}
\begin{array}{ll}
{b}^T \mathcal{A}_{\ell,t}^{\hat{{x}}}({x}) {b} \leq \lambda_{\ell} {b}^T M_{\ell,t}^{\hat{{x}}} {b}, & t = 1,\ldots,n_k, \quad \mathrm{for~all~} {b} \in \mathbb{R}^{N_{\ell}},\\[1ex]
{c}^T \mathcal{A}_{u,t}^{\hat{{x}}}({x}) {c} \geq \lambda_u {c}^T M_{u,t}^{\hat{{x}}} {c}, & t = 1,\ldots,n_k, \quad \mathrm{for~all~} {c} \in \mathbb{R}^{N_u}.\\[1ex]
\end{array}
\end{equation}
We will approximate the above conditions  by judiciously generating a finite number of approximating vectors $b^{(1)}, \ldots, b^{(N_B)}  \in \mathbb{R}^{N_{\ell}}$ and $c^{(1)}, \ldots, c^{(N_C)} \in \mathbb{R}^{N_{u}}$.  (The method for choosing and updating these sets of vectors will be discussed in the next subsection.)  The resulting linear inequalities in the variables $x$, $\lambda_{\ell}$, and $\lambda_u$ are:
\begin{equation}\label{eq_sdp2lp2}
\begin{array}{ll}
(b^{(p)})^T A_{\ell,t,0}(b^{(p)}) + \sum_i^{n_{\varepsilon}} (b^{(p)})^T A_{\ell,t,i}(b^{(p)}) x_i \leq \lambda_{\ell} (b^{(p)})^T M_{\ell,t}^{\hat{{x}}} (b^{(p)}), & t = 1,\ldots,n_k, \quad p=1, \ldots, N_B \ ,\\[1ex]
(c^{(p)})^T A_{u,t,0}(c^{(p)}) + \sum_i^{n_{\varepsilon}} (c^{(p)})^T A_{u,t,i}(c^{(p)}) x_i \ge \lambda_{u} (c^{(p)})^T M_{u,t}^{\hat{{x}}} (c^{(p)}), & t = 1,\ldots,n_k, \quad p=1, \ldots, N_C \ .\\[1ex]
\end{array}
\end{equation}

Because the mass matrices $M_{\ell,t}^{\hat{{x}}}$ and $M_{u,t}^{\hat{{x}}}$ are positive definite (Proposition \ref{theorem_RedMat}), the coefficients in the right-hand-side of \eqref{eq_sdp2lp2} are all positive.  It follows that \eqref{eq_sdp2lp2} can be reformatted by rescaling as the following two linear inequality systems:
\begin{equation}\label{eqLinRelax}
B^{\hat{{x}}}{x} + g^{\hat{{x}}} \le {e} \lambda_{\ell}\ , \qquad C^{\hat{{x}}}{x} + h^{\hat{{x}}} \ge {e} \lambda_u \ ,
\end{equation}
where $B^{\hat{{x}}}\in\mathbb{R}^{(N_B n_k)\times n_x}$, and $C^{\hat{{x}}}\in\mathbb{R}^{(N_C n_k)\times n_x}$.  It then follows from Proposition \ref{theorem_RedMat} that:
\begin{equation}\label{prop_positive}
(B^{\hat{x}} x + g^{\hat{x}})_j > 0, \quad j = 1, \ldots, N_B n_k =: \mathcal{N}_B \ \ \mathrm{and} \ \ (C^{\hat{x}} x + h^{\hat{x}})_i > 0, \quad i = 1, \ldots, N_C n_k =: \mathcal{N}_C \ ,
\end{equation} for all $x \in S$.


Replacing the semidefinite inclusions in \eqref{eqSDP} with their linear inequality approximations \eqref{eqLinRelax}, we obtain the following linear fractional approximation of \eqref{eqSDP}:
\begin{equation}\label{eq_lfpapp}
\begin{array}{lclr}
P_{LFP}^{\hat x}: &  \underset{ x,\lambda_{\ell},\lambda_u}{\max} & \displaystyle \frac{\lambda_u - \lambda_{\ell}}{\lambda_u+\lambda_{\ell}}& \\ \\
     & \mbox{ s.t. } &  B^{\hat{{x}}}{x} + g^{\hat{{x}}}\le {e}\lambda_{\ell}, & \\ [1.5ex]
     &               &  C^{\hat{{x}}}{x} + h^{\hat{{x}}}\ge {e} \lambda_u, & \\ [1.5ex]
     &               &  x_{\min} \leq x_i \leq x_{\max} \ , i=1, \ldots, n_x\\[1.5ex]
     &  &  \lambda_{\ell} \ge 0, \ \lambda_u \ge 0, \ \lambda_{\ell} + \lambda_u > 0 \ .
\end{array}
\end{equation}
The superscript ``$(\cdot)^{\hat{{x}}}$'' indicates that components of $B$, $C$, $g$, and $h$ are functions of (and so depend on) ${\hat{{x}}}$. In order for the linear inequality formulation to be reasonably accurate, the optimal solution ${x}^{\ast}$ of \eqref{eq_lfpapp} should be close enough to the linearizing point ${\hat{{x}}}$, i.e., $\|{x}^{\ast} - {\hat{{x}}}\|\le \epsilon$.  Since the optimization problem \eqref{eq_lfpapp} is a linearly constrained linear fractional optimization problem, it can be converted to a linear program and efficiently solved by using standard linear optimization software. Table \ref{tb_LP} presents the basic outline of the algorithm for solving bandgap optimization problems by the linear fractional optimization \eqref{eq_lfpapp} instead of the semidefinite program \eqref{eqSDP}.  We note  in Step 4 of the algorithm  that one can augment the solution process for $P^{\hat x}_{LFP}$ with a standard delayed constraint generation procedure \citep{bt}.  More detailed implementation of Step 2b is discussed in the next subsection.

\begin{table}
\caption{Algorithm for solving bandgap problems using linear inequalities approximation of eigenvalue bounds.\label{tb_LP}}
\begin{center}
\begin{tabular}{|ll|}
\hline
& \textbf{Algorithm for Bandgap Optimization using Linear Inequalities Approximation}\\\hline
{\bf Step 1.} & Start with initial guess $\hat x := x^0$ and tolerance $\epsilon_{\mbox{tol}}$\\  & \\
{\bf Step 2a.}& Construct the matrices \eqref{eqRedMat} for \eqref{eqSDP} based on $\hat x$, for  each $t\in\{1, \ldots, n_k\}$: \\
              & \quad ${\cal A}_{\ell,t}^{\hat{{x}}}(x) := A_{\ell,t,0}^{\hat{x}} + \sum_{i=1}^{n_x} A_{\ell,t,i}^{\hat{x}} x_i$\\
              & \quad ${\cal A}_{u,t}^{\hat{{x}}}( x) := A_{u,t,0}^{\hat{x}} + \sum_{i=1}^{n_x} A_{u,t,i}^{\hat{x}} x_i $\\
              & \quad $M_{\ell,t}^{\hat{{x}}}$ and $M_{u,t}^{\hat{{x}}}$\\ & \\
{\bf Step 2b.}& Choose vectors $b^1, \ldots, b^{N_B}$ and $c^1, \ldots, c^{N_C}$:\\  &\\
{\bf Step 2c.}& Construct the data for \eqref{eq_lfpapp} based on $\hat x$ and the linear operators from Step 2a:\\
              & \quad $B^{\hat{{x}}}$, $g^{\hat{x}}$, $C^{\hat{{x}}}$, and $h^{\hat{x}}$\\ & \\
{\bf Step 3.} & Form the linear fractional problem $P^{\hat x}_{LFP}$ in \eqref{eq_lfpapp} \\ & \\
{\bf Step 4.} & Solve $P^{\hat x}_{LFP}$ for an optimal solution $(x^*, \lambda_{\ell}^*, \lambda_u^*)$\\
              & \quad {\em (Optional:  augment $P^{\hat x}_{LFP}$ with Delayed Constraint Generation}) \\ & \\
{\bf Step 5.} & If $\|x^* - \hat x\| \le \epsilon_{\mbox{tol}}$, stop.\\
              & Else update $\hat x \leftarrow x^*$ and go to {\bf Step 2.}\\
\hline
\end{tabular}
\end{center}
\end{table}

\subsection{Methodology for Constructing the Approximating Vectors}
\label{subsec_vec}

We describe our approach for constructing the approximating vectors $b^{(1)}, \ldots, b^{(N_B)}  \in \mathbb{R}^{N_{\ell}}$ and $c^{(1)}, \ldots, c^{(N_C)} \in \mathbb{R}^{N_{u}}$.  We focus on $b^{(1)}, \ldots, b^{(N_B)}  \in \mathbb{R}^{N_{\ell}}$, as the same approach is also used to construct the approximating vectors $c^{(1)}, \ldots, c^{(N_C)} \in \mathbb{R}^{N_{u}}$.  Note that ${N_{\ell}}$ (and $N_u$) is not large, typically ${N_{\ell}} \approx 3-7$, due to the subspace approximation. Ideally, we would want the approximating vectors to be distributed uniformly over the upper half of the Euclidean sphere: $\{ b \in \mathbb{R}^{N_{\ell}} : \sqrt{b^Tb} = 1, \ b_{N_{\ell}} \ge 0\}$, where we need only consider a half-sphere because $v^TMv = (-v)^TM(-v)$ for any $v \in \mathbb{R}^{N_{\ell}}$.  For ease of construction, we choose to work with the upper half of the unit $L_1$-sphere, also known as the upper boundary of the cross-polytope $\{ b \in \mathbb{R}^{N_{\ell}} : \|b\|_1 = 1, \ b_{N_{\ell}} \ge 0\}$, whose $2^{(N_{\ell}-1)}$ facets are the unit $(N_{\ell}-1)$-simplices in their respective orthants.  Given an integer dilation constant $K$, we first define:
\begin{equation}\label{eq_K}
{\cal K} := \Bigl\{ {k} \in \mathbb{R}^{N_{\ell}} : \sum_{i=1}^{N_{\ell}} |{k}_i| = K, \ {k}_i \text{ integer } \Bigr\} \ ,
\end{equation}
and then define the approximating vectors $b^{(1)}, \ldots, b^{(N_B)}  \in \mathbb{R}^{N_{\ell}}$ to be the elements of the following set:
$$ {B}_K :=\Bigl\{ {b}\in \mathbb{R}^{N_{\ell}}: {b} =(1/K){k} \text{ for some } {k} \in {\cal K},\ {k}_{N_{\ell}} \ge 0 \ \Bigr\}.$$
The resulting approximating vectors are distributed uniformly on the surface of the half cross-polytope. This is illustrated in Figure \ref{fig_basisVector} for $N_{\ell} = 2$. Note that  the number of vectors in $ {B}_K$ grows as $ O(K^{N_{\ell}-1})$. Increasing $K$ will render the piecewise linear approximation model more accurate albeit at higher computational cost. In addition and if necessary, we expand the set of approximating vectors at each iteration using delayed constraint generation: once the linear fractional optimization problem \eqref{eq_lfpapp} is solved, we check the semidefinite inclusions in \eqref{eqSDP}  for any eigenvectors violating the constraints and add them to the set of approximating vectors to generate additional linear inequality cuts which are then added to \eqref{eq_lfpapp}.  Note that checking the semidefinite inclusions in \eqref{eqSDP} is inexpensive due to the reduced size of the system.
\begin{figure}
\centering
\includegraphics[scale=0.65]{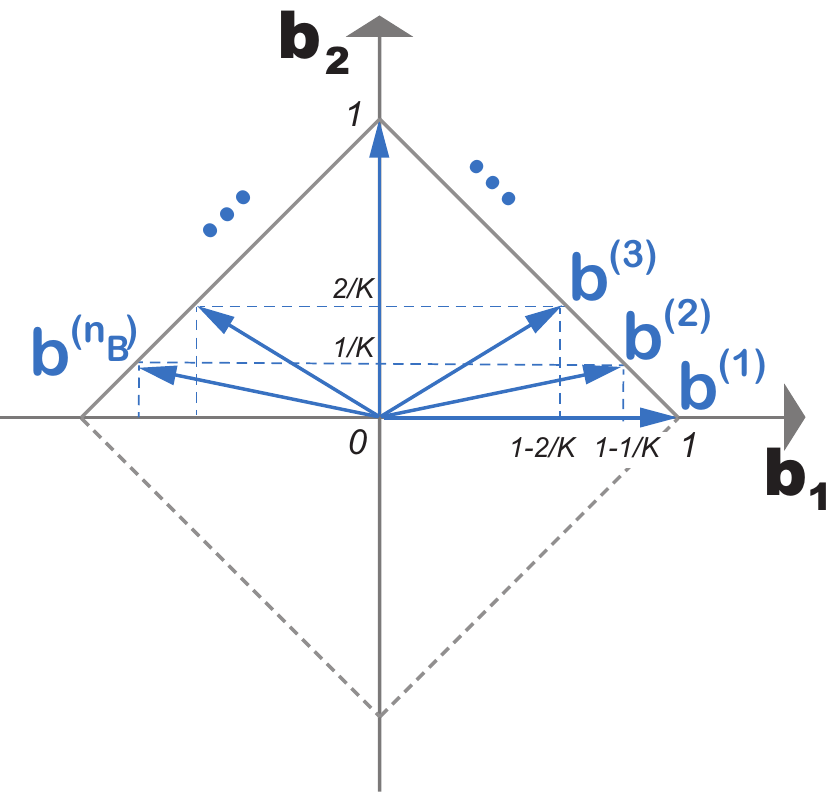}
\caption{The vectors chosen to construct the approximating linear inequalities are distributed uniformly on the surface of a half cross-polytope arising from the $L_1$ norm.}
\label{fig_basisVector}
\end{figure}


\subsection{Quality of Linear Inequalities Approximation}\label{subsec_finality}

To verify the quality of the approximation of the semidefinite inclusions by using the above approach, we focus on the effect of the tunable parameter $K$ (defined in equation \eqref{eq_K}) and the resulting number of linear inequalities. We first note that when $K$ is large, more vectors (larger $N_B$ and $N_C$) are generated to presumably approximate semidefinite inclusions  more accurately, yielding a linear fractional optimization problem \eqref{eq_lfpapp} that better approximates the semidefinite problem \eqref{eqSDP}.  As a result, the resulting linear optimization problem will contain a larger number of linear inequalities and thus require more computation time. On the other hand, a smaller value of $K$ will reduce the number of linear inequality constraints, but result in a less accurate approximation of \eqref{eqSDP}.  The quality of the linear inequalities approximation may be empirically measured in terms of the number of outer iterations of the algorithm of Table \ref{tb_LP} and the number of ``successful'' solutions, where a solution is deemed successful if it opens up a bandgap more than $10\%$.\medskip

We conduct an empirical test in order to determine a good value of $K$. In particular, we make $10$ runs of the algorithm of Table \ref{tb_LP} using $10$ randomly chosen starting point configurations for a variety of types of bandgap problems and report the results in Table \ref{tab_outerIter}.  In this table, the headings in the right side columns of the form $\Delta^{TE}_{1,2}$ refer to bandgap optimization of the bandgap between the $2$nd and $1$st eigenvalues in TE polarization, etc.  Table \ref{tab_outerIter} shows average outer iterations, and number of successful runs, for various bandgap optimization problems by using the algorithm of Table \ref{tb_LP}, with a large value of $K$ ($K=5$, resulting in $N_B, N_C \sim 500$) and a small value of $K$ ($K=3$, resulting in $N_B, N_C \sim 10 $) combined with delayed constraint generation (DCG). The results of the SDP approach (presented in \cite{men2010bandgap}) are shown in the table as a benchmark for comparison.  We observe that using a small value of $K=3$ combined with delayed constraint generation appears to strike a good compromise between system size (and computation time) and the success rate.

\begin{table}
\caption{Average number of outer iterations and the total number of successful runs (out of $10$ runs) for bandgap optimization using the semidefinite program (SDP) formulation and the linear fractional program (LFP) formulation. $\Delta\lambda_{m,m+1}$ indicates the optimized bandgap is between the $m^{\mathrm{th}}$ and the $(m+1)^{\mathrm{st}}$ eigenvalues. DCG denotes delayed constraint generation.}
\label{tab_outerIter}
\centering
\begin{subtable}{\linewidth}\centering
\renewcommand{\arraystretch}{1.5}
{\begin{tabular}{|l||c|c|c|c|c|}
\hline
{Bandgap} & $\Delta\lambda_{1,2}^{TE}$  & $\Delta\lambda_{2,3}^{TE}$ & $\Delta\lambda_{8,9}^{TE}$ & $\Delta\lambda_{9,10}^{TE}$\\
\hline\hline
{\textbf{\em{SDP}}}  & $9.0/\emph{7}$ & $9.0/\emph{6}$ & $14.2/ \emph{2}$ &$23.5/ \emph{1}$\\
\hline
{\textbf{\em{LFP}}}  ($K = 5$) & $17.0/\emph{9}$ & $9.1/\emph{8}$ & $43.5/\emph{3}$ & $40.5/\emph{3}$ \\
{\textbf{\em{LFP}}} ($K = 3$) & $20.0/\emph{6}$ & $12.6/\emph{6}$ & $37.1/\emph{1}$ & $26.1/\emph{2}$  \\
{\textbf{\em{LFP}}} ($K = 3$) with DCG & $14.0/\emph{8}$  & $15.8/\emph{6}$ & $27.7/\emph{3}$ & $24.1/\emph{4}$ \\
\hline
\end{tabular}}
\label{tab_outerIter_TE}
\caption{TE polarization}
\end{subtable}
\begin{subtable}{\linewidth}\centering
\renewcommand{\arraystretch}{1.5}
{\begin{tabular}{|l||c|c|c|c|c|}
\hline
{Bandgap} &  $\Delta\lambda_{1,2}^{TM}$ & $\Delta\lambda_{2,3}^{TM}$ & $\Delta\lambda_{8,9}^{TM}$ & $\Delta\lambda_{9,10}^{TM}$\\
\hline\hline
{\textbf{\em{SDP }}}  & $3.4/\emph{10}$ & $4.1/\emph{8}$ & $10.9/\emph{3}$ & $22.5/\emph{2}$\\
\hline
{\textbf{\em{LFP}}}  ($K = 5$) & $5.1/\emph{10}$ & $10.2/\emph{7}$ & $31.2/\emph{3}$ & $36.1/\emph{4}$ \\
{\textbf{\em{LFP}}}  ($K = 3$) & $5.2/\emph{10}$ & $6.3/\emph{8}$ & $20.5/\emph{2}$ & $34.2/\emph{2}$ \\
{\textbf{\em{LFP}}}  ($K = 3$) with DCG & $5.2/\emph{10}$ & $6.9/\emph{7}$ & $23.2/\emph{2}$ & $27.6/\emph{2}$ \\
\hline
\end{tabular}}
\label{tab_outerIter_TM}
\caption{TM polarization}
\end{subtable}
\end{table}

 \end{APPENDICES}

\ACKNOWLEDGMENT{This work is supported by AFOSR Grant No. FA9550-11-1-0141, the Singapore-MIT Alliance, the MIT-Chile-Pontificia Universidad Católica de Chile Seed Fund, and LaCaixa Fellowship.}


\bibliographystyle{ormsv080}
\bibliography{fr_ref}






%
%

\end{document}